\documentclass[10pt]{amsart}
\usepackage[margin=2.5cm]{geometry}
\usepackage{tikz-cd}
\usepackage{graphicx}
\usepackage{mathdots}		
\usepackage{amssymb}
\usepackage{amsmath}
\usepackage{relsize}
\usepackage{subfig, graphicx}
\usepackage{mathrsfs}
\usepackage{amssymb, amsthm, indentfirst}
\usepackage{relsize}
\usepackage{hyperref}
\usepackage{stmaryrd}
\usepackage{lineno}
\usepackage{mathabx}
\numberwithin{equation}{section}
\usepackage{color}
\theoremstyle{plain}
\newtheorem{thm}{Theorem}[section]

\newtheorem{conj}[thm]{Conjecture}
\newtheorem{lemma}[thm]{Lemma}
\newtheorem{prop}[thm]{Proposition}
\newtheorem{corollary}[thm]{Corollary}

\theoremstyle{definition}

\newtheorem*{example*}{Example}
\newtheorem{rmk}[thm]{Remark}
\def\Gal{\operatorname{Gal}}

\newtheorem*{hypothesis*}{Hypothesis}

\makeatletter
\DeclareFontFamily{OMX}{MnSymbolE}{}
\DeclareSymbolFont{MnLargeSymbols}{OMX}{MnSymbolE}{m}{n}
\SetSymbolFont{MnLargeSymbols}{bold}{OMX}{MnSymbolE}{b}{n}
\DeclareFontShape{OMX}{MnSymbolE}{m}{n}{
    <-6>  MnSymbolE5
   <6-7>  MnSymbolE6
   <7-8>  MnSymbolE7
   <8-9>  MnSymbolE8
   <9-10> MnSymbolE9
  <10-12> MnSymbolE10
  <12->   MnSymbolE12
}{}
\DeclareFontShape{OMX}{MnSymbolE}{b}{n}{
    <-6>  MnSymbolE-Bold5
   <6-7>  MnSymbolE-Bold6
   <7-8>  MnSymbolE-Bold7
   <8-9>  MnSymbolE-Bold8
   <9-10> MnSymbolE-Bold9
  <10-12> MnSymbolE-Bold10
  <12->   MnSymbolE-Bold12
}{}

\let\llangle\@undefined
\let\rrangle\@undefined
\DeclareMathDelimiter{\llangle}{\mathopen}%
                     {MnLargeSymbols}{'164}{MnLargeSymbols}{'164}
\DeclareMathDelimiter{\rrangle}{\mathclose}%
                     {MnLargeSymbols}{'171}{MnLargeSymbols}{'171}
\makeatother

\author{Shih-Yu Chen}
\address{Department of Mathematics, National Tsing Hua University, 101, Section 2, Kuang-Fu Road, Hsinchu 300, Taiwan, ROC}
\email{sychen.math@gmail.com}

\def\SL{{\rm{SL}}}
\def\GL{{\rm{GL}}}
\def\GSp{{\rm GSp}}
\def\GU{{\rm GU}}
\def\U{{\rm U}}

\def\A{{\mathbb A}}
\def\C{{\mathbb C}}

\def\K{{\mathcal K}}

\def\R{{\mathbb R}}
\def\Q{{\mathbb Q}}
\def\Z{{\mathbb Z}}

\def\<{\langle}
\def\>{\rangle}

\def\bp{\begin{pmatrix}}
\def\ep{\end{pmatrix}}

\def\<{\langle}
\def\>{\rangle}

\def\GL{\operatorname{GL}}

\def\GSp{\operatorname{GSp}}

\def\SL{\operatorname{SL}}

\def\U{\operatorname{U}}

\def\1{\mathbf{1}}

\def\itPi{\mathit{\Pi}}
\def\itPsi{\mathit{\Psi}}
\def\itSigma{\mathit{\Sigma}}

\makeatletter
\newcommand{\extp}{\@ifnextchar^\@extp{\@extp^{\,}}}
\def\@extp^#1{\mathop{\bigwedge\nolimits^{\!#1}}}
\makeatother

\makeatletter
\newcommand{\exterior}[1]{\mathop{\mathpalette\exterior@{#1}}}
\newcommand{\exterior@}[2]{%
  \raisebox{\depth}{%
  \fontsize{\sf@size}{0}%
  \m@th
  $\ifx#1\displaystyle\textstyle\else#1\fi\bigwedge$}%
  ^{\mspace{-2mu}#2}%
  \kern-\scriptspace
}
\makeatother

\title{Algebraicity of adjoint $L$-functions for quasi-split groups}
\begin{document}

\begin{abstract}
For a globally generic cuspidal automorphic representation $\itPi$ of a quasi-split reductive group $G$ over $\Q$, E. Lapid and Z. Mao proposed a conjecture on the decomposition of global Whittaker functionals on $\itPi$ into products of an adjoint $L$-value of $\itPi$ and local Whittaker functionals. 
In this paper, we consider the algebraic aspect of the Lapid--Mao conjecture. More precisely, when $\itPi$ is $C$-algebraic, we show that the algebraicity of the adjoint $L$-value can be expressed in terms of the Petersson norm of Whittaker-rational cusp forms in $\itPi$, subject to the validity of the Lapid--Mao conjecture.
For unitary similitude groups, we also establish an unconditional and more refined algebraicity result. 
Additionally, we give an explicit formula for the case $G={\rm U}(2,1)$.
\end{abstract}

%% \tableofcontents %% Just for papers exceeding 50 pages.
\maketitle
%\pagewiselinenumbers
\tableofcontents

\section{Introduction}

Let $G$ be a quasi-split connected reductive group over $\Q$. Let $N$ be the unipotent radical of a Borel subgroup of $G$ over $\Q$. Let $\A$ be the ring of rational adeles and fix a non-degenerated character $\psi_N$ of $N(\Q)\backslash N(\A)$. 
For an automorphic form $f$ on $G(\A)$, consider its Whittaker function with respect to $\psi_N$ defined by
\[
W_{\psi_N}(g;f) = \int_{N(\Q)\backslash N(\A)}f(ng)\overline{\psi_{N}(n)}\,dn^{\rm Tam},\quad g \in G(\A)
\]
where $dn^{\rm Tam}$ is the Tamagawa measure.
Let $\itPi$ be a globally $\psi_N$-generic cuspidal automorphic representation of $G(\A)$. 
Recall that globally $\psi_N$-generic means the Whittaker functional $W_{\psi_N}$ does not vanishing on $\itPi$. In this case, $\itPi$ admits a unique model, called the Whittaker model of $\itPi$, realized in the space of Whittaker functions on $G(\A)$ with respect to $\psi_N$ which is denoted by $\mathcal{W}(\itPi,\psi_N)$. 
For a sufficiently large finite set $S$ of places of $\Q$ containing the archimedean place, let
\[
L^S(s,\itPi,{\rm Ad}) = \prod_{v \notin S}L^S(s,\itPi_v,{\rm Ad}),\quad {\rm Re}(s) \gg 0
\]
be the partial adjoint $L$-function of $\itPi$, and $\Delta_G(s)$ the $L$-function of the dual Gross motive associated to $G$ (cf.\,\cite{Gross1997}).
It is expected that the ratio
\[
\frac{L^S(s,\itPi,{\rm Ad})}{\Delta_G^S(s)}
\]
admits meromorphic continuation to the whole complex plane and is holomorphic and non-vanishing at $s=1$.
In \cite{LM2015}, Lapid and Mao proposed a Ichino--Ikeda type formula (cf.\,{\S\,\ref{SS:LM} below) which relates the value of the ratio at $s=1$ to the global Whittaker period. Suitably interpret, it expresses the ratio at $s=1$ in terms of the Petersson pairing
\begin{align}\label{E:Petersson norm}
\<f,f'\> = \int_{G(\Q)\backslash G(\A)^1}f(g)f'(g)\,dg,\quad f \in \itPi,\,f' \in \itPi^\vee
\end{align}
where $G(\A)^1$ is the intersection of the kernel of $|\chi|_\A$ for all $\chi \in X^*(G)$, and the Haar measure $dg$ is normalized so that ${\rm vol}(G(\Q)\backslash G(\A)^1,dg)=1$.
In $loc.$ $cit.$, it is verified for $G = {\rm Res}_{F/\Q}\GL_n$ (see also \cite{LM2017} when $G$ is a metaplectic group). When $G = {\rm Res}_{F/\Q}\GSp_4$, it is proved by Furusawa and Morimoto in \cite[Theorem 6.3]{FM2024}. Recently, by Beuzart-Plessis and Chaudouard \cite{BPC2023} and Morimoto \cite{Morimoto2024}, the formula is completely proved when $G$ is (the Weil restriction to 
$\Q$ of) a unitary group (see also \cite{Morimoto2022}).

In this paper, we consider the algebraic aspect of the Lapid--Mao conjecture.
More precisely, recall that the rationality field $\Q(\itPi_f)$ of $\itPi_f := \otimes_p' \itPi_p$ is the fixed field in $\C$ of the set of automorphisms $\sigma \in {\rm Aut}(\C)$ such that ${}^\sigma\!\itPi_f \cong \itPi_f$. 
Similar notation applies to $\itPi^S := \otimes_{p \notin S }\itPi_p$ and $\itPi_p$ for any prime number $p$.
In \cite{BG2014}, Buzzard and Gee introduced notions of been ``algebraic'' for automorphic representations of $G(\A)$, which generalize the corresponding notion for the case $G = {\rm Res}_{F/\Q}\GL_n$ introduced by Clozel \cite{Clozel1990}. 
More precisely, we say $\itPi$ is $C$-algebraic if the infinitesimal character of $\itPi_\infty$, after twisted by the Weyl vector, is integral. 
For instance, $\itPi$ is $C$-algebraic if it is cohomological.
If $\itPi$ is $C$-algebraic, then it is conjectured in \cite[Conjecture 5.15]{BG2014} that $\itPi^S$ is defined over a number field for sufficiently large $S$, which implies that $\Q(\itPi^S)$ is a number field. In fact, we expect that in this case $\itPi_f$ would be defined over a number field (holds if $\itPi$ is cohomological, see \cite[Proposition 2.4]{GR2013}).
Now we assume that
\begin{align}\label{E:assumption}
\mbox{$\itPi$ is $C$-algebraic and $\itPi_\infty$ is essentially tempered}.
\end{align}
Here essentially tempered means being tempered after twisting by a character in {$X^*(G)\otimes_\Z\Q$}. 
Let $\itPi_+$ be the space of cusp forms in $\itPi$ whose archimedean component are highest weight vectors of the minimal $K_\infty$-type of $\itPi_\infty$, where $K_\infty$ is a maximal compact subgroup of $G(\R)$. 
For a subfield $E$ of $\C$ containing $\Q(\itPi_f)\cdot\Q^{\rm ab}$, we have the notion of been Whittaker-rational over $E$ for cusp forms in $\itPi_+$ (with respect to a choice of highest weight Whittaker function of $\itPi_\infty$ in the minimal $K_\infty$-type, see \S\,\ref{SS:algebraicity}).
The contragredient representation $\itPi^\vee$ is globally $\psi_N^{-1}$-generic and also satisfies assumption (\ref{E:assumption}).
Let $\itPi_-^\vee$ be the space of cusp forms in $\itPi^\vee$ whose archimedean component are lowest weight vectors of the minimal $K_\infty$-type of $\itPi_\infty^\vee$. Similarly we have the notion of been Whittaker-rational for cusp forms in $\itPi_-^\vee$. 
Subject to the validity of the Lapid--Mao conjecture, we prove the following result on the algebraicity of adjoint $L$-values of $\itPi$ in terms of Petersson pairing of Whittaker-rationals cusp forms.

\begin{thm}[Theorem \ref{T: main 0}]\label{T:main}
Assume (\ref{E:assumption}) and the Lapid--Mao conjecture hold for $\itPi$.
Let $E$ be a subfield of $\C$ containing $\Q(\itPi_f)\cdot \Q^{\rm ab}$. 
For $f \in \itPi_+$ and $f' \in \itPi^\vee_-$ with $\<f,f'\>\neq 0$ which are Whittaker-rational over $E$, we have
\[
\left.\left(\frac{L^S(s,\itPi,{\rm Ad})}{\Delta_G^S(s)}\right)\right\vert_{s=1}\frac{1}{C_{N,\itPi_\infty}\cdot\<f,f'\>} \in E.
\]
Here $C_{N,\itPi_\infty}$ is a non-zero complex number depending on $N$ and $\itPi_\infty$ defined in (\ref{E:C_infty 1}).
\end{thm}

\begin{rmk}
We can replace (\ref{E:assumption}) by the assumption that $\itPi_\infty$ admits a unique minimal $K_\infty$-type. Nonetheless, the above result is ``algebraic'' only when $E/\Q$ is an algebraic extension in which case $\itPi$ is expected to be $C$-algebraic.
The essentially temperedness of $\itPi_\infty$ is part of the generalized Ramanujan conjecture suggested by Howe and Piatetski-Shapiro.
\end{rmk}

The following natural questions arise from the theorem:
\begin{itemize}
    \item[(1)] Descend the algebraicity from $E$ to a number field, for instance to $\mathbb{Q}(\itPi_f)$ itself.
    \item[(2)] Compute explicitly the archimedean constant $C_{N,\itPi_\infty}$.
\end{itemize}

For question (1), additional global input is usually required to guarantee the existence of a model over $\mathbb{Q}(\itPi_f)$. This is known, for example, when either $\itPi$ is cohomological or $G$ admits a Shimura datum. The ambiguity arising from $\mathbb{Q}^{\mathrm{ab}}$ is due to the potential non-uniqueness of the Whittaker datum.
Question (2) is purely local in nature, but the computation of $C_{N,\itPi_\infty}$ is notoriously difficult in general.
Several special cases have been treated in the literature, as we now recall.
When $G = \mathrm{Res}_{F/\mathbb{Q}}\mathrm{GL}_2$ and $\itPi$ is cohomological, both questions have affirmative answers (cf.\,\cite[Proposition 5.15]{Namikawa2015} and the references therein). For higher-rank general linear groups over number fields, question (1) was also known by Balasubramanyam and Raghuram in \cite{BR2017} if $\itPi$ is cohomological. In this setting, the local constant $C_{N,\itPi_\infty}$ was computed in \cite[Lemmas 5.3 and 5.4]{Chen2020} and \cite[Lemma 6.4]{Chen2021} for $\mathrm{GL}_3(\mathbb{R})$ and $\mathrm{GL}_3(\mathbb{C})$, respectively.
When $G = \mathrm{Res}_{F/\mathbb{Q}}\mathrm{GSp}_4$ for a totally real field $F$ and $\itPi_\infty$ is a (limit of) discrete series representation, question (1) was answered in our earlier work \cite{Chen2021c}, and the local problem (2) was approached via global methods with Ichino in \cite{CI2019}.
In the above-mentioned cases, having an explicit formula (or at least one valid modulo the constant $C_{N,\itPi_\infty}$) plays a crucial role in relating the adjoint $L$-value to the congruence module.
For instance, Hida’s seminal work on congruences between modular forms \cite{Hida1981} and its generalizations to $\mathrm{GL}_2$ over number fields by Urban \cite{Urban1995}, Ghate \cite{Ghate2002}, Dimitrov \cite{Dimitrov2005}, and Namikawa \cite{Namikawa2015};
the result of Balasubramanyam--Raghuram \cite{BR2017} for $\mathrm{GL}_n$ over number fields and cohomological $\itPi$;
the result of Lemma--Ochiai \cite{LO2018} for $\mathrm{GSp}_4$ and endoscopic $\itPi$;
and the recent work of Ricoul \cite{Ricoul2025} for $\mathrm{GL}_3$ over real quadratic fields and conjugate self-dual (or self conjugate) $\itPi$ all exemplify this philosophy.

%Besides the potential application of Theorem \ref{T:main} to congruences between automorphic representations, 
Theorem \ref{T:main}, together with its refinement obtained after resolving questions (1) and (2), also has applications to the study of the algebraicity of special $L$-values.
More precisely, when $G$ admits a Shimura datum and $\itPi_\infty$ is a (limit of) discrete-series representation, it is expected that the Petersson pairing of Whittaker-rational cusp forms decomposes as a product of two automorphic periods arising from coherent cohomology.
In a pervious work \cite{Chen2021f}, this expectation is worked out in detail for $G = \mathrm{GSp}_4$, and the resulting period relation is applied to prove Deligne’s conjecture for the symmetric $6$-th $L$-functions of modular forms, as well as the algebraicity of the critical values of $L$-functions for $\mathrm{GSp}_4 \times \mathrm{GSp}_4$.
In a recent work of Horawa and Prasanna \cite{HP2025}, the period relation in the limit-of-discrete-series case is used in their study of the motivic action on the Siegel modular threefold.
In Theorem \ref{T: main 1} below, we resolve Question (1) for unitary similitude groups when $\itPi_\infty$ is a (limit of) discrete series representation.
In Proposition \ref{P:archimedean}, we perform the archimedean computation for $\mathrm{U}(2,1)$, and an explicit formula is given in Theorem \ref{T:explicit}.
These results will be applied to study the motivic action on the Picard modular surface in an ongoing joint project with Horawa.

\section{Algebraicity of adjoint $L$-values}

\subsection{Lapid--Mao conjecture}\label{SS:LM}

In this section, we recall the conjecture proposed by Lapid and Mao in \cite{LM2015}. 
Let $G$ be a quasi-split connected reductive group over $\Q$. Let $N$ be the unipotent radical of a Borel subgroup of $G$ over $\Q$, and fix a non-degenerated character $\psi_N$ of $N(\Q)\backslash N(\A)$. 
For each prime number $p$, we fix a speical maximal compact subgroup $K_p$ of $G(\Q_p)$ which is hyperspecial if $G$ is unramified at $p$. Let $dn_p$ be the Haar measure on $N(\Q_p)$ normalized so that ${\rm vol}(N(\Q_p)\cap K_p,dn_p)=1$.
For the archimedean place $\infty$, let $dn_\infty$ be the Haar measure on $N(\R)$ such that $\prod_v dn_v$ is the Tamagawa measure on $N(\A)$.

Let $\itPi$ be a globally $\psi_N$-generic cuspidal automorphic representation of $G(\A)$.
Consider the global Whittaker functional $I$ in 
$
{\rm Hom}_{N(\A)\times N(\A)}(\itPi\otimes\itPi^\vee,\C)
$
defined by
\[
I(f\otimes f') = W_{\psi_N}(1;f)W_{\psi_N^{-1}}(1;f').
\]
For each place $v$ of $\Q$, we fix a $G(\Q_v)$-equivariant bilinear pairing $\<\cdot,\cdot\>_v$ between $\itPi_v$ and $\itPi^\vee_v$.
When $v=p$ is non-archimedean and $\itPi_p$ is unramified, we fix non-zero $K_p$-invariant vectors $f_{\itPi_p} \in \itPi_p$ and $f_{\itPi_p^\vee} \in \itPi_p^\vee$.
We assume that the pairings $\<\cdot,\cdot\>_v$ are chosen so that for $f = \otimes_vf_v \in \itPi$ and $f' = \otimes_vf_v' \in \itPi^\vee$ (with respect to the fixed unramified vectors), we have $\<f_v,f_v'\>_v=1$ for almost all $v$ and
\[
\<f,f'\> = \prod_v \<f_v,f_v'\>_v.
\]
Let $I_v \in {\rm Hom}_{N(\Q_v)\times N(\Q_v)}(\itPi_v\otimes\itPi_v^\vee,\C)$ be the local Whittaker functional defined by
the stable integral
\begin{align}\label{E:local Whittaker 1}
I_v(f_v\otimes f_v') = \int_{N(\Q_v)}^{\rm st}\<\itPi_v(n_v)f_v,f_v'\>_v\overline{\psi_{N,v}(n_v)}\,dn_v
\end{align}
when $v$ is non-archimedean, and we refer to \cite[\S\,2.5]{LM2015} for the definition when $v$ is archimedean.
For instance, if $\itPi_\infty$ is a discrete series representation, then $I_\infty$ is defined as above with stable integral replaced by the usual integration. 
%We define the normalized functional $I_v^*$ by
%\[
%I_v^* = \frac{L(1,\itPi_v,{\rm Ad})}{L(1,\omega_{\K_v/\Q_v})\zeta_v(2)L(3,\omega_{\K_v/\Q_v})}\cdot I_v.
%\] 
If $v=p$ is non-archimedean, $\itPi_p$ and $\psi_{N,p}$ are unramified, then we have (cf.\,\cite[Proposition 2.14]{LM2015})
\begin{align}\label{E:unramified LM}
\frac{I_p(f_{\itPi_p}\otimes f_{\itPi_p^\vee})}{\<f_{\itPi_p}, f_{\itPi_p^\vee}\>_p} = \frac{\Delta_{G,p}(1)}{L(1,\itPi_p,{\rm Ad})}.
\end{align}

\begin{conj}[Lapid--Mao]\label{C:LM}
As elements of ${\rm Hom}_{N(\A)\times N(\A)}(\itPi\otimes\itPi^\vee,\C)$, we have
\begin{align}\label{E:LM}
I = \frac{1}{|\mathcal{S}_{\itPi}|}\cdot \prod_v I_v,
\end{align}
where $\mathcal{S}_\itPi$ is the global component group of the (conjectural) Arthur parameter of $\itPi$.
\end{conj}

\subsection{Algebraicity of adjoint $L$-values}\label{SS:algebraicity}

We keep the notation of \S\,\ref{SS:LM}.
For each place $v$ of $\Q$, we denote by $\mathcal{W}(\itPi_v,\psi_{N,v})$ the Whittaker model of $\itPi_v$ with respect to $\psi_{N,v}$ which consisting of smooth functions $W$ on $G(\Q_v)$ such that 
\[
W(ng) = \psi_{N,v}(n)W(g),\quad g \in G(\Q_v),\,n \in N(\Q_v)
\]
and are of moderate growth if $v$ is archimedean.
When $v=p$ is non-archimedean, $\itPi_p$ and $\psi_{N,p}$ are unramified, let $W_{\itPi_p} \in \mathcal{W}(\itPi_p,\psi_{N,p})$ be the right $K_p$-invariant Whittaker function normalized so that $W_{\itPi_p}(1) = 1$. Note that the value at $1$ of the unramified Whittaker functions are non-zero by the result of Casselman and Shalika \cite[Theorem 5.4]{CS1980}. Similarly, we denote by $\mathcal{W}(\itPi_f,\psi_{N,f})$ the Whittaker model of $\itPi_f$ with respect to $\psi_{N,f}$. Note that we have a $G(\A_f)$-equivariant isomorphism
\[
\otimes_p'\mathcal{W}(\itPi_p,\psi_{N,p}) \longrightarrow \mathcal{W}(\itPi_f,\psi_{N,f}),\quad \otimes_p W_p \longmapsto \prod_p W_p.
\]
Assume $\pi$ is $C$-algebraic and $\itPi_\infty$ is essentially tempered. As we recalled in the introduction, in this case we expect that $\itPi_f$ is defined over a number field. 
{Since $\itPi_\infty$ is essentially tempered with real infinitesimal character (such representations are called \textit{temperic} after A. Afgoustidis), by the result of Vogan \cite[Theorem 1.2]{Vogan2007}, $\itPi_\infty$ admits a unique minimal $K_\infty$-type for any maximal compact subgroup $K_\infty$ of $G(\R)$.}
By replacing $K_\infty$ with a conjugate of it if necessary, there exists a highest weight Whittaker function $W_{\infty,+} \in \mathcal{W}(\itPi_\infty,\psi_{N,\infty})$ in the minimal $K_\infty$-type such that
\[
W_{\infty,+}(1) \neq 0.
\]
Such Whittaker functions are then unique up to scalars. 
Let $\itPi_+$ be the space of cusp forms in $\itPi$ whose archimedean component are highest weight vectors of the minimal $K_\infty$-type of $\itPi_\infty$.
With respect to a fixed choice of $W_{\infty,+}$, we have the $G(\A_f)$-equivariant isomorphism
\[
\mathcal{W}(\itPi_f,\psi_{N,f}) \longrightarrow \itPi_+,\quad W \longmapsto f_W
\]
defined so that $f_W$ is the unique cusp form in $\itPi_+$ such that
\begin{align}\label{E:Whittaker identification}
W_{\psi_N}(f_W) = W_{\infty,+}\cdot W.
\end{align}
Let $\Q^{\rm ab}$ be the maximal abelian extension of $\Q$ in $\C$. For $\sigma \in {\rm Aut}(\C/\Q(\itPi_f)\cdot\Q^{\rm ab})$, we define 
\[
{}^\sigma W(g):=\sigma(W(g)),\quad g \in G(\A_f).
\]
Since $\psi_{N,f}$ takes values in $\Q^{\rm ab}$ and ${}^\sigma\!\itPi \cong \itPi$, we have ${}^\sigma W \in \mathcal{W}(\itPi_f,\psi_{N,f})$.
We thus have a $\sigma$-linear $G(\A_f)$-equivariant isomorphism 
\[
\itPi_+ \longrightarrow \itPi_+,\quad f\longmapsto {}^\sigma\!f
\]
defined by ${}^\sigma\!(f_W):=f_{{}^\sigma W}$. For a subfield $E$ of $\C$ containing $\Q(\itPi_f)\cdot\Q^{\rm ab}$, we say $f \in \itPi_+$ is Whittaker-rational over $E$ if ${}^\sigma\!f = f$ for all $\sigma \in {\rm Aut}(\C/E)$.
The contragredient representation $\itPi^\vee$ is globally $\psi_N^{-1}$-generic and also satisfies assumption (\ref{E:assumption}).
Let $\itPi_-^\vee$ be the space of cusp forms in $\itPi^\vee$ whose archimedean component are lowest weight vectors of the minimal $K_\infty$-type of $\itPi_\infty^\vee$. Fix a lowest weight Whittaker function $W_{\infty,-} \in \mathcal{W}(\itPi_\infty^\vee,\psi_{N,\infty}^{-1})$ in the minimal $K_\infty$-type such that $W_{\infty,-}(1) \neq 0$. Similarly we have the notion of been Whittaker-rational for cusp forms in $\itPi_-^\vee$ (with respect to the choice $W_{\infty,-}$). 
We define an archimedean constant
\begin{align}\label{E:C_infty 1}
C_{N,\itPi_\infty}:= W_{\infty,+}(1)^{-1}W_{\infty,-}(1)^{-1}\cdot\frac{I_\infty(f_{\infty,+}\otimes f_{\infty,-})}{\<f_{\infty,+},f_{\infty,-}\>_\infty},
\end{align}
where $f_{\infty,+} \in \itPi_\infty$ (resp.\,$f_{\infty,-} \in \itPi_\infty^\vee$) is a highest (resp.\,lowest) weight vector in the minimal $K_\infty$-type. 
It is clear that the ratio in the right-hand side of (\ref{E:C_infty 1}) is independent of the choices of $f_{\infty,+}$, $f_{\infty,-}$, and $\<\cdot,\cdot\>_\infty$. 
The local dependence on $\itPi_\infty$ is due to the choices of Whittaker functions $W_{\infty,+}$ and $W_{\infty,-}$.
The global dependence on $N$ is due to the choice of Haar measure $dn_\infty$ so that $\prod_v dn_v$ is the Tamagawa measure.

\begin{thm}\label{T: main 0}
Assume $\itPi$ is $C$-algebraic, $\itPi_\infty$ is essentially tempered, and Conjecture \ref{C:LM} holds for $\itPi$.  
For $f \in \itPi_+$ and $f' \in \itPi^\vee_-$ with $\<f,f'\>\neq 0$, we have
\[
\sigma \left(\left.\left(\frac{L^S(s,\itPi,{\rm Ad})}{\Delta_G^S(s)}\right)\right\vert_{s=1}\frac{ 1}{C_{N,\itPi_\infty}\cdot\<f,f'\>} \right) = \left.\left(\frac{L^S(s,\itPi,{\rm Ad})}{\Delta_G^S(s)}\right)\right\vert_{s=1}\frac{1}{C_{N,\itPi_\infty}\cdot\<{}^\sigma\!f,{}^\sigma\!f'\>}
\]
for all $\sigma \in {\rm Aut}(\C/\Q(\itPi_f)\cdot\Q^{\rm ab})$, where $S$ is a sufficiently large finite set of places of containing the archimedean place.
In particular, if $f$ and $f'$ are Whittaker-rational over a subfield $E$ of $\C$ containing $\Q(\itPi_f)\cdot\Q^{\rm ab}$, then we have
\[
\left.\left(\frac{L^S(s,\itPi,{\rm Ad})}{\Delta_G^S(s)}\right)\right\vert_{s=1}\frac{1}{C_{N,\itPi_\infty}\cdot\<f,f'\>} \in E.
\]
\end{thm}

\subsection{Proof of Theorem \ref{T: main 0}}\label{SS: main 0 pf}

In this section, we prove Theorem \ref{T: main 0} subject to the validity of the Lapid--Mao conjecture.
Let $p$ be a prime number and $\sigma \in {\rm Aut}(\C/\Q^{\rm ab})$. We fix $\sigma$-linear $G(\Q_p)$-equivariant isomorphisms
\begin{align*}
\itPi_p &\longrightarrow {}^\sigma\!\itPi_p,\quad f_p \longmapsto {}^\sigma\!f_p\\
\itPi_p^\vee &\longrightarrow {}^\sigma\!\itPi_p^\vee,\quad f_p' \longmapsto {}^\sigma\!f_p'.
\end{align*} 
Let $\<\cdot,\cdot\>_{p,\sigma}$ be the $G(\Q_p)$-equivariant pairing between ${}^\sigma\!\itPi_p$ and ${}^\sigma\!\itPi_p^\vee$ defined by
\[
\<{}^\sigma\!f_p,{}^\sigma\!f_p'\>_{p,\sigma}:=\sigma\left(\<f_p,f_p'\>_p\right),\quad f_p \in \itPi_p,\,f_p'\in\itPi_p^\vee.
\]
We then define the local Whittaker functional $I_{p,\sigma} \in {\rm Hom}_{N(\Q_p)\times N(\Q_p)}({}^\sigma\!\itPi_p\otimes{}^\sigma\!\itPi_p^\vee,\C)$ as in (\ref{E:local Whittaker 1}) using the pairing $\<\cdot,\cdot\>_{p,\sigma}$. 
The following lemma is on the Galois equivariant property of the local Whittaker functional.

\begin{lemma}\label{L:Galois equiv. LM 0}
Let $p$ be a prime number and $\sigma \in {\rm Aut}(\C/\Q^{\rm ab})$. 
For $f_p \in \itPi_p$ and $f_p' \in \itPi_p^\vee$, we have
\[
\sigma\left(I_p(f_p\otimes f_p')\right) = I_{p,\sigma}({}^\sigma\!f_p \otimes {}^\sigma\!f_p').
\]
\end{lemma}

\begin{proof}
By assumption, ${}^\sigma\!\psi_{N,p} = \psi_{N,p}$.
In particular, for any open compact subgroup $U_p$ of $N(\Q_p)$
\begin{align*}
&\sigma \left ( \int_{U_p} \<\itPi_p(n_p)f_p,f_p'\>_p\overline{\psi_{N,p}(n_p)}\,dn_p\right)\\
& = \int_{U_p}\sigma\left(\<\itPi_p(n_p)f_p,f_p'\>_p\right)\overline{\sigma(\psi_{N,p}(n_p))}\,dn_p\\
& = \int_{U_p}\<{}^\sigma\!\itPi_p(n_p){}^\sigma\!f_p,{}^\sigma\!f_p'\>_{p,\sigma}\overline{\psi_{N,p}(n_p)}\,dn_p.\end{align*}
The assertion then follows immediately from the definition of stable integral.
\end{proof}

%For $\sigma \in {\rm Aut}(\C/\Q(\itPi_f)\cdot \Q^{\rm ab})$, the map
%\[
%\mathcal{W}(\itPi_f,\psi_{N,f}) \longrightarrow \mathcal{W}(\itPi_f,\psi_{N,f}),\quad W \longmapsto {}^\sigma W
%\]
%is a $\sigma$-linear $G(\A_f)$-equivariant isomorphism.
%For any subfield $E$ of $\C$ containing $\Q(\itPi_f)$, we denote by $\mathcal{W}(\itPi_f,\psi_{N,f})^{{\rm Aut}(\C/E\cdot \Q^{\rm ab})}$ the space of ${\rm Aut}(\C/E\cdot \Q^{\rm ab})$-invariant Whittaker functions.
%By definition, $f_W \in \itPi_+$ is Whittaker-rational over $E\cdot \Q^{\rm ab}$ if and only if $W \in \mathcal{W}(\itPi_f,\psi_{N,f})^{{\rm Aut}(\C/E\cdot \Q^{\rm ab})}$.

%\begin{lemma}
%Let $E$ be a field of definition of $\itPi_f$. Then $\mathcal{W}(\itPi_f,\psi_{N,f})^{{\rm Aut}(\C/E\cdot \Q^{\rm ab})}$ defines a $E\cdot \Q^{\rm ab}$-rational structure of $\mathcal{W}(\itPi_f,\psi_{N,f})$ if and only if it is non-zero.
%\end{lemma}

%\begin{proof}

%\end{proof}

Now we begin the proof of Theorem \ref{T: main 0}.
Let $S$ be a finite set of places containing the archimedean place such that $\itPi_p$ and $\psi_{N,p}$ are both unramified for $p \notin S$. Put $T = S \smallsetminus \{\infty\}$ and we write
\[
\itPi_T = \otimes_{p \in T}\itPi_p,\quad \psi_{N,T} = \prod_{p \in T}\psi_{N,p},\quad \<\cdot,\cdot\>_T = \prod_{p \in T}\<\cdot,\cdot\>_{p},\quad I_T = \prod_{p \in T}I_p.
\]
To prove Theorem \ref{T: main 0}, it suffices to assume
\[
f = f_{\infty,+} \otimes f_T \otimes (\otimes_{p \notin T}f_{\itPi_p}),\quad f' = f_{\infty,-} \otimes f_T' \otimes (\otimes_{p \notin T}f_{\itPi_p^\vee})
\]
for some $f_T \in \itPi_T$ and $f_T' \in\itPi_T^\vee$ such that $\<f_T,f_T'\>_T\neq0$.
For $p \notin T$, let $W_{\itPi_p} \in \mathcal{W}(\itPi_p,\psi_{N,p})$ and $W_{\itPi_p^\vee} \in \mathcal{W}(\itPi_p,\psi_{N,p}^{-1})$ be the $K_p$-invariant Whittaker functions normalized so that $W_{\itPi_p}(1) = W_{\itPi_p^\vee}(1)=1$.
Fix $W_T \in \mathcal{W}(\itPi_T,\psi_{N,T})$ and $W_T' \in \mathcal{W}(\itPi_T^\vee,\psi_{N,T}^{-1})$  so that $W_T(1)$ and $W_T'(1)$ are non-zero. Let
\[
W = W_T \cdot \prod_{p\notin T} W_{\itPi_p} \in \mathcal{W}(\itPi_f,\psi_{N,f}),\quad W' = W_T' \cdot \prod_{p\notin T} W_{\itPi_p^\vee}\in \mathcal{W}(\itPi_f^\vee,\psi_{N,f}^{-1}).
\]
Then we have
\[
f_W = f_{\infty,+} \otimes g_T \otimes (\otimes_{p \notin T}f_{\itPi_p}),\quad f_{W'} = f_{\infty,-} \otimes g_T' \otimes (\otimes_{p \notin T}f_{\itPi_p^\vee})
\]
for some $g_T \in \itPi_T$ and $g_T'\in\itPi_T^\vee$.
By (\ref{E:unramified LM}) and (\ref{E:LM}) applying to $I(f_W\otimes f_{W'})$, we have
\begin{align}\label{E:LM main pf 1}
\frac{W_T(1)W_T'(1)}{\<f,f'\>} = \frac{1}{|\mathcal{S}_\itPi|}\cdot \left.\left(\frac{\Delta_G^S(s)}{L^S(s,\itPi,{\rm Ad})}\right)\right\vert_{s=1}\cdot C_{N,\itPi_\infty} \cdot \frac{I_T(g_T\otimes g_T')}{\<f_T,f_T'\>_T}.
\end{align}
Let $\sigma \in {\rm Aut}(\C/\Q(\itPi_f)\cdot \Q^{\rm ab})$.
Note that
\[
{}^\sigma\!f = f_{\infty,+} \otimes {}^\sigma\!f_T \otimes (\otimes_{p \notin T}f_{\itPi_p}),\quad {}^\sigma\!f' = f_{\infty,-}\otimes {}^\sigma\!f_T' \otimes (\otimes_{p \notin T}f_{\itPi_p^\vee}).
\]
Similarly for ${}^\sigma\!(f_W)$ and ${}^\sigma\!(f_{W'})$.
Apply (\ref{E:LM}) to $I({}^\sigma\!(f_W)\otimes {}^\sigma\!(f_{W'}))$, we have
\begin{align}\label{E:LM main pf 2}
\frac{{}^\sigma W_T(1){}^\sigma W_T'(1)}{\<{}^\sigma\!f,{}^\sigma\!f'\>} = \frac{1}{|\mathcal{S}_\itPi|}\cdot \left.\left(\frac{\Delta_G^S(s)}{L^S(s,\itPi,{\rm Ad})}\right)\right\vert_{s=1}\cdot C_{N,\itPi_\infty} \cdot \frac{I_{T,\sigma}({}^\sigma\!g_T\otimes {}^\sigma\!g_T')}{\<{}^\sigma\!f_T,{}^\sigma\!f_T'\>_{T,\sigma}}.
\end{align}
By Lemma \ref{L:Galois equiv. LM 0}, 
\[
\sigma \left( \frac{I_T(g_T\otimes g_T')}{\<f_T,f_T'\>_T} \right) = \frac{I_{T,\sigma}({}^\sigma\!g_T\otimes {}^\sigma\!g_T')}{\<{}^\sigma\!f_T,{}^\sigma\!f_T'\>_{T,\sigma}}.
\]
The assertion then follows from applying $\sigma$ to both sides of (\ref{E:LM main pf 1}) and then compare with (\ref{E:LM main pf 2}).
This completes the proof of Theorem \ref{T: main 0}.

\section{Adjoint $L$-values for unitary similitude groups}

\subsection{Groups, Whittaker data}\label{SS:basic setting}
Let $\K$ be an imaginary quadratic field of discriminant $-D_\K$. We write $x\mapsto x^c$ for the non-trivial automorphism of $\K$.
Let $\iota$ and $\overline{\iota} = \iota\circ c$ be the embeddings of $\K$ into $\C$. We fix a CM-type $\{\iota\}$ and identify $\K$ with a subfield of $\C$ via $\iota$.
For $n \geq 1$, let $Q_n \in \GL_n$ defined by
\[
Q_n = \begin{cases}
\bp  & {\bf 1}_{n/2} \\ -{\bf 1}_{n/2} &  \ep & \mbox{ if $n$ is even},\\
\bp & & 1 \\ &\iddots& \\ 1 & &  \ep  & \mbox{ if $n$ is odd}.
\end{cases}
\]
Let ${\rm GU}_n$ be the quasi-split unitary similitude group over $\Q$ defined by
\[
{\rm GU}_n = \left\{ g \in {\rm Res}_{\K/\Q}\GL_n \,\left\vert\, {}^t\!g^c Q_n g = \nu(g)Q_n,\,\nu(g) \in \mathbb{G}_m   \right.\right\}.
\]
Let ${\rm U}_n$ be the kernel of the similitude character $\nu$. It is well-known that ${\rm GU}_n$ admits a unique Whittaker datum up to conjugation. We fix a choice of Whittaker datum $(B_n,\psi_{N_n})$ as follows.
Let $B_n$ be the Borel subgroup of ${\rm GU}_n$ defined by 
\[
B_n = \left.\left\{ \bp a & \\ & {}^ta^{-c} \ep \bp {\bf 1}_{n/2} & b \\ & {\bf 1}_{n/2}\ep \,\right\vert\, \mbox{$a$ is upper triangular in ${\rm Res}_{\K/\Q}\GL_{n/2}$},\, b={}^tb^c  \right\}
\]
if $n$ is even, and consisting of upper triangular matrices in $\GU_n$ if $n$ is odd. 
Let $N_n$ be the unipotent radical of $B_n$, and $\psi_{N_n}$ be the non-degenerated character of $N_n(\Q)\backslash N_n(\A)$ defined by
\[
\psi_{N_n}(u) = \begin{cases}
\psi\circ {\rm Tr}_{\K/\Q}\left(\sum_{i=1}^{n/2-1}u_{i,i+1} + \tfrac{1}{2}u_{n/2,n}\right) & \mbox{ if $n$ is even},\\
\psi\circ {\rm Tr}_{\K/\Q}\left(\sum_{i=1}^{(n-1)/2}u_{i,i+1} \right) & \mbox{ if $n$ is odd},
\end{cases}
\]
where $\psi$ is the non-trivial additive character of $\Q\backslash\A$ such that $\psi_\infty(x)=e^{2\pi\sqrt{-1}\,x}$.

\subsection{$C$-algebraic discrete series representations}\label{SS:d.s. rep.}

For $r,s \in\Z_{\geq 0}$ with $r+s>0$, let $\GU(r,s)$ be the real unitary similitude group of signature $(r,s)$ defined by
\[
\GU(r,s) = \left\{ g \in \GL_{r+s}(\C) \,\left\vert\, {}^t\overline{g}\bp {\bf 1}_r & \\ & -{\bf 1}_s \ep g = \nu(g)\bp {\bf 1}_r & \\ & -{\bf 1}_s \ep,\,\nu(g) \in \R^\times   \right.\right\}.
\]
We denote by $\GU(r,s)^+$ the subgroup consisting of elements with positive similitude.
With respect to the fixed CM-type $\{\iota\}$, we identify $\K_\infty = \K\otimes_\R\C$ with $\C$. Note that the archimedean signature of $\GU_n$ is $(\tfrac{n+\delta}{2},\tfrac{n-\delta}{2})$, where $\delta \in \{0,1\}$ with $\delta \equiv n \,({\rm mod }\,2)$. Fix $g_\infty \in \GL_n(\C)$ such that 
\begin{align}\label{E:polarization}
 (\sqrt{-1})^{1-\delta}\cdot Q_n = {}^t\overline{g}_\infty \bp {\bf 1}_r & \\ & -{\bf 1}_s\ep g_\infty.
\end{align}
From now on we identify $\GU_n(\R)$ with $\GU(\tfrac{n+\delta}{2},\tfrac{n-\delta}{2})$ via the isomorphism
\[
\GU_n(\R) \longrightarrow \GU(\tfrac{n+\delta}{2},\tfrac{n-\delta}{2}),\quad g \longmapsto g_\infty g g_\infty^{-1}.
\]
Let $K_\infty$ be the maximal connected modulo center compact subgroup of $\GU_n(\R)$ consisting of block diagonal matrices ${\rm diga}(ak_1,ak_2)$ with $k_1 \in \U(\tfrac{n+\delta}{2})$, $k_2 \in \U(\tfrac{n-\delta}{2})$, and $a \in \R^\times$.
Let $H_\infty$ be the maximal torus of $\GU_n(\R)$ consisting of diagonal matrices. Then $H_\infty = T_\infty\R^\times$, where $T_\infty \cong \U(1)^n$ is the maximal compact subgroup of $H_\infty$. We denote by $\frak{g}$, $\frak{k}$, and $\frak{h}$ the Lie algebras of $\GU_n(\R)$, $K_\infty$, and $H_\infty$ respectively, and $\frak{g}_\C$, $\frak{k}_\C$, and $\frak{h}_\C$ be their complexifications. 
We identify $\frak{h}_\C^*$ with $\C^{n}\times\C$ by the map
\[
z_1e_1+\cdots +z_ne_n + wf \longmapsto (z_1,\cdots,z_n;w),
\]
where
\[
e_i({\rm diag}(t_1,\cdots,t_n)) = \begin{cases}
t_i & \mbox{ if $t_j \in \sqrt{-1}\,\R$ for all j},\\
0 & \mbox{ if $t_1= \cdots = t_n \in \R$},
\end{cases}
\]
and 
\[
f({\rm diag}(t_1,\cdots,t_n)) = \begin{cases}
0 & \mbox{ if $t_j \in \sqrt{-1}\,\R$ for all j},\\
t_1 & \mbox{ if $t_1= \cdots = t_n \in \R$}.
\end{cases}
\]
The set of roots of $\frak{g}_\C$ with respect to $\frak{h}_\C$ is then given by
\[
R = \left\{ e_i-e_j \,\vert\, 1 \leq i \neq j \leq n \right\}.
\]
We denote by $R_c$ (resp.\,$R_{nc}$) the set of compact (resp.\,non-compact) roots.
For $\tau \in \{\iota,\overline{\iota}\}$, let $R_\tau^+$ be a positive system of $(\frak{g}_\C,\frak{h}_\C)$ defined by $R_{\overline{\iota}}^+ = -R_{\iota}^+$ and
\[
R_\iota^+ = \left\{ e_i-e_j \,\vert\, 1 \leq i < j \leq n \right\}.
\]
Let $\delta_\tau$ be half the sum of roots in $R_\tau^+$. It is clear that $\delta_\iota = -\delta_{\overline{\iota}}$.
Let $R_{\tau,c}^+ = R_c \cap R_\tau^+$ (resp.\,$R_{\tau,nc}^+ = R_{nc} \cap R_\tau^+$) be the set of positive compact (resp.\,non-compact) roots with respect to $\tau$. 
We denote by $\frak{p}_\tau^+$ and $\frak{p}_\tau^-$ the spaces of $\frak{g}_\C$ generated by the root spaces of roots in $R_{\tau,nc}^+$ and $-R_{\tau,nc}^+$ respectively. Let $\frak{P}_\tau$ be the parabolic subalgebra of $\frak{g}_\C$ defined by $\frak{P}_\tau = \frak{k}_\C \oplus \frak{p}_\tau^-$, and $P_\tau$ the closed subgroup of $\GU_n(\C)$ with Lie algebra $\frak{P}_\tau$.

For $\Lambda = (\lambda;{\sf w}) \in \C^n\times\C$ with $\lambda = (\lambda_1,\cdots,\lambda_n)$, define $\Lambda^\vee = (\lambda^\vee;-{\sf w})$ and $\Lambda^c = (-\lambda;{\sf w})$ with
\[
\lambda_i^\vee = \begin{cases}
-\lambda_{(n+\delta)/2-i+1} & \mbox{ if $1 \leq i \leq \tfrac{n+\delta}{2}$},\\
-\lambda_{n+(n+\delta)/2-i+1} & \mbox{ if $\tfrac{n+\delta}{2} < i \leq n$}.
\end{cases}
\]
When $n$ is even, we also define $\tilde{\Lambda} = (\tilde{\lambda};{\sf w})$ by interchanging the first and last $n/2$ entries of $\lambda$:
\[
\tilde{\lambda} = (\lambda_{n/2+1},\cdots,\lambda_n,\lambda_1,\cdots,\lambda_{n/2}).
\]
We say $\Lambda \in \C^n \times \C$ is an admissible weight if it is $R$-regular and $\Lambda = (\lambda;{\sf w}) \in (\Z+\tfrac{n-1}{2})^n\times\Z$ for some $\lambda$ and ${\sf w}$ such that 
\begin{align}\label{E:d.s. parity}
\lambda_1+\cdots+\lambda_n \equiv {\sf w}\,({\rm mod}\,2).
\end{align}
Recall that $C$-algebraic discrete series representations of $\GU_n(\R)^+$ are parameterized by the $W(\frak{k}_\C,\frak{h}_\C)$-orbits of admissible weights.
More precisely, for an admissible weight $\Lambda = (\lambda;{\sf w})$, we have a discrete series representation
$\pi_\lambda$ of $\U_n(\R)$ constructed by Harish-Chandra (cf.\,\cite[IX, \S\,7]{Knapp1986}).
The parity condition (\ref{E:d.s. parity}) implies that we can extend the representation to $\pi_\Lambda$ on $\GU_n(\R)^+$ by 
\[
\pi_\Lambda(ag) := a^{\sf w}\pi_\lambda(g),\quad a \in \R^\times,\,g\in\U_n(\R).
\]
We call $\Lambda$ the Harish-Chandra parameter of $\pi_\Lambda$.
Moreover, $\pi_\Lambda \cong \pi_{\Lambda'}$ if and only if $\Lambda$ and $\Lambda'$ are conjugate by $W(\frak{k}_\C,\frak{h}_\C)$.
Note that we have
\[
\pi_\Lambda^\vee \cong \pi_{\Lambda^\vee},\quad \pi_\Lambda^c \cong \pi_{\Lambda^c}.
\]
Also note that if $n$ is even, then $\pi_\Lambda$ is isomorphic to the conjugate of $\pi_{\tilde{\Lambda}}$ by an element in $\GU_n(\R)$ of negative similitude.
%It is clear that $\pi_\Lambda$ is $C$-algebraic if and only if $\Lambda = (\lambda;{\sf w})$ for some ${\sf w}\in\Z$.
Let $\tau \in \{\iota,\overline{\iota}\}$ and consider the $(\frak{P}_\tau,K_\infty)$-cohomology of discrete series representations (cf.\,\cite[Theorem 4.6.2]{Harris1990} and \cite[Theorem 3.2.1]{BHR1994}). An admissible weight $\Lambda$ uniquely determines a positive system $R_{\Lambda}^+$ such that $\Lambda$ is positive with respect to it. Assume further that $\Lambda$ is $R_{\tau,c}^+$-positive, then the following $(\frak{P}_\tau,K_\infty)$-cohomology group is one-dimensional:
\[
H^{q(\Lambda)}(\frak{P}_\tau,K_\infty; \pi_{\Lambda} \otimes V_{\Lambda^\vee - \delta_\tau}),
\]
where $q(\Lambda) = |R_\Lambda^+\cap(-R^+_{\tau,nc})|$ and $V_{\Lambda^\vee - \delta_\tau}$ is the irreducible algebraic representation of $K_\infty$ of $R_{\tau,c}^+$-highest weight $\Lambda^\vee - \delta_\tau$.
Note that $\Lambda^\vee$ is also $R_{\tau,c}^+$-positive.
Moreover, if $\Lambda'$ is another $R_{\tau,c}^+$-positive admissible weight such that the $(\frak{P}_\tau,K_\infty)$-cohomology of $\pi_{\Lambda'} \otimes V_{\Lambda^\vee - \delta_\tau}$ is non-vanishing at some degree $q$, then we have $\Lambda' = \Lambda$ and $q=q(\Lambda)$.
Here we only concern with $C$-algebraic generic discrete series representations, which are of the form $\pi_\Lambda$ for admissible $R_{\iota,c}^+$-positive weight $\Lambda$ such that the simple roots of $R_{\Lambda}^+$ are all non-compact. Quantitatively, this condition is equivalent to $\Lambda = (\lambda;{\sf w})$ with 
\begin{align}\label{E:generic parameter}
\lambda = (\alpha_1,\alpha_3,\cdots, \alpha_{n-1+\delta}; \alpha_2,\alpha_4,\cdots,\alpha_{n-\delta}),\quad \alpha_1>\alpha_2>\cdots>\alpha_n
\end{align}
or if $n$ is even we can also have
\[
\lambda = (\alpha_2,\alpha_4,\cdots,\alpha_{n};\alpha_1,\alpha_3,\cdots, \alpha_{n-1}),\quad \alpha_1>\alpha_2>\cdots>\alpha_n.
\]
When $n$ is odd, we have $\GU_n(\R) = \GU_n(\R)^+$. In this case $C$-algebraic discrete series representations of $\GU_n(\R)$ are constructed as above, and we write $\itPi_\Lambda = \pi_\Lambda$.
When $n$ is even, for an admissible weight $\Lambda$, we denote by $\itPi_\Lambda$ the unique $C$-algebraic discrete series representation of $\GU_n(\R)$ such that its restriction to $\GU_n(\R)^+$ contains $\pi_\Lambda$. Then we have
\[
\itPi_\Lambda\vert_{\GU_n(\R)^+} \cong \pi_\Lambda \oplus \pi_{\tilde{\Lambda}}.
\]
Thus in this case $\itPi_\Lambda \cong \itPi_{\tilde{\Lambda}}$.

\subsection{Coherent cohomology of unitary Shimura varieties}

Let $({\rm GU}_n,X_\tau)$ be the Shimura datum with $X_\tau$ equal to the $\GU_n(\R)$-conjugacy class containing the morphism $h_\tau : {\rm Res}_{\C/\R}\mathbb{G}_m \rightarrow \GU_{n,\R}$ with
\[
h_\tau(z) = \begin{cases}
{\rm diag}(z{\bf 1}_{\small{\frac{n+\delta}{2}}}, \overline{z}{\bf 1}_{\scriptstyle\frac{n-\delta}{2}}) & \mbox{ if $\tau = \iota$},\\
{\rm diag}(\overline{z}{\bf 1}_{\scriptstyle\frac{n+\delta}{2}}, {z}{\bf 1}_{\scriptstyle\frac{n-\delta}{2}}) & \mbox{ if $\tau = \overline{\iota}$},
\end{cases}
\]
on $\R$-points. It is clear that $K_\infty$ is the stabilizer of $h_\tau$.
We have the associated unitary Shimura variety 
\[
{\rm Sh}^\tau = \varprojlim_{K} {\rm Sh}^\tau_{K} = \varprojlim_{{K}}\GU_n(\Q)\backslash X_\tau\times \GU_n(\A_{f}) / {K},
\]
where ${K}$ runs through neat open compact subgroups of $\GU_n(\A_{f})$.
It is a pro-algebraic variety over $\C$ with continuous $\GU_n(\A_{f})$-action and admits a canonical model over the reflex field $E({\rm GU}_n,X_\tau)$. Note that the reflex field is equal to $\Q$ (resp.\,$\K$) if $n$ is even (resp.\,odd).

Fix $\tau \in \{\iota,\overline{\iota}\}$. Let $(\rho,V_\rho)$ be an irreducible algebraic representation of $K_\infty$. We extend it to a representation of $P_\tau$ so that it factors through the reductive quotient $K_{\infty,\C}$ of $P_\tau$.
The resulting representation then determines a $\GU_{n,\C}$-vector bundle $\mathcal{V}_\rho$ over the flag variety $\GU_n(\C) / P_\tau$ (cf.\,\cite[(1.3.1)]{BHR1994}).
Let $\beta_\tau^*(\mathcal{V}_\rho)$ be the pullback $\GU_n(\R)$-bundle over $X_\tau$, where $\beta_\tau$ is the Borel embedding
\[
\beta_\tau : X_\tau \cong \GU_n(\R) / K_\infty \longrightarrow \GU_n(\C)/P_\tau.
\]
We then obtain a vector bundle 
\[
[\mathcal{V}_\rho] = \varprojlim_{{K}}\,[\mathcal{V}_\rho]_K = \varprojlim_{{K}} \GU_n(\Q)\backslash \,\beta_\tau^*(\mathcal{V}_\rho) \times \GU_n(\A_f) / {K}
\]
over the Shimura variety ${\rm Sh}^\tau$. In the literature (cf.\,\cite{Harris1985}, \cite{Harris1986}, \cite{Milne1990}), this is called the automorphic vector bundle over ${\rm Sh}^\tau$ associated with $\rho$. Note that the complexified representation $\rho_\C$ of $K_{\infty,\C}$ is always defined over $\K$. In particular, by the results of Harris \cite[Theorem 3.3]{Harris1985} and Milne \cite[Theorem 5.1]{Milne1990}, $[\mathcal{V}_\rho]$ admits a canonical model over $\K$. 
For a neat open compact subgroup $K$ and a good admissible smooth rational polyhedral
cone decomposition data $\Sigma$, we have the toroidal compactification ${\rm Sh}_{K,\Sigma}^\tau$ of ${\rm Sh}_K^\tau$. 
We remark that for the PEL-type Shimura varieties such as ${\rm Sh}^\tau$, we have a rather algebraic construction of the toroidal compactifications by Lan \cite{Lan2013}.
By the results of Harris \cite[Proposition 2.8]{Harris1989} and Pink \cite[\S\,12.4]{Pink1989}, ${\rm Sh}_{K,\Sigma}^\tau$ admits a canonical model over the reflex field. Moreover, by \cite[\S\,4]{Harris1989}, $[\mathcal{V}_\rho]_K$ admits a canonical extension $[\mathcal{V}_\rho]_{K,\Sigma}^{\rm can}$ over ${\rm Sh}_{K,\Sigma}^\tau$ which is defined over $\K$.
The subcanonical extension $[\mathcal{V}_\rho]_{K,\Sigma}^{\rm sub}$ is define to be the tensor of $[\mathcal{V}_\rho]_{K,\Sigma}^{\rm can}$ with the ideal sheaf defining the divisor ${\rm Sh}_{K,\Sigma}^\tau \smallsetminus {\rm Sh}_{K}^\tau$.
Let
\begin{align*}
H^\bullet({\rm Sh}^\tau,[\mathcal{V}_\rho]^{\rm can})&:= \varinjlim_{K,\Sigma} H^\bullet({\rm Sh}_{K,\Sigma}^\tau,[\mathcal{V}_\rho]_{K,\Sigma}^{\rm can}),\\
H^\bullet({\rm Sh}^\tau,[\mathcal{V}_\rho]^{\rm sub})&:= \varinjlim_{K,\Sigma} H^\bullet({\rm Sh}_{K,\Sigma}^\tau,[\mathcal{V}_\rho]_{K,\Sigma}^{\rm sub})
\end{align*}
be the direct limit of coherent cohomology of $[\mathcal{V}_\rho]_{K,\Sigma}^{\rm can}$ or $[\mathcal{V}_\rho]_{K,\Sigma}^{\rm sub}$ over ${\rm Sh}_{K,\Sigma}^\tau$, where the $\Sigma$'s are directed by refinement.
We denote by $H^\bullet_!({\rm Sh}^\tau,[\mathcal{V}_\rho])$ (interior cohomology) the image of $H^\bullet({\rm Sh}^\tau,[\mathcal{V}_\rho]^{\rm sub})$ in $H^\bullet({\rm Sh}^\tau,[\mathcal{V}_\rho]^{\rm can})$.
%These cohomology groups natural $\GU_n(\A_f)$-action inherited from that on the automorphic vector bundles.
%Also the action is admissible (cf.\,\cite[Proposition 2.6]{Harris1990}).
The cohomology groups naturally carry an action of $\GU_n(\A_f)$, inherited from that on the automorphic vector bundles. Moreover, this action is admissible, as shown in \cite[Proposition 2.6]{Harris1990}. 
By considering the conjugates of Shimura varieties and automorphic vector bundles under ${\rm Aut}(\C)$ as explained in \cite[Proposition 1.4.3 and \S\,4]{BHR1994}, for each $\sigma \in {\rm Aut}(\C)$, we have a $\sigma$-linear $\GU_n(\A_f)$-equivariant isomorphism
\begin{align*}
T_\sigma : H^\bullet_!({\rm Sh}^\tau,[\mathcal{V}_\rho]) \longrightarrow H^\bullet_!({\rm Sh}^{\sigma\circ\tau},[\mathcal{V}_{{}^\sigma\!\rho}]),
\end{align*}
where 
\[
{}^\sigma\!\rho = \begin{cases}
\rho & \mbox{ if $\sigma \vert_\K = {\rm id}$},\\
\rho^c & \mbox{ if $\sigma \vert_\K = c$}.
\end{cases}
\]
%We will recall in the proof of Proposition \ref{P:rational structure} some relations of the coherent cohomology with the $(\frak{P}_\tau,K_\infty)$-chomology of automorphic forms on $\GU_n(\A)$.
%Let $\tau \in \{\iota,\overline{\iota}\}$ and $V_\rho$ an irreducible algebraic representation of $K_\infty$.
In \cite[Corollary 3.4]{Harris1990}, Harris proved that the coherent cohomology group $H^\bullet({\rm Sh}^\tau,[\mathcal{V}_\rho]^{\rm can})$ can be computed by the $(\frak{P}_\tau,K_\infty)$-cohomology of smooth functions on $\GU_n(\Q)\backslash \GU_n(\A)$ of moderate growth with coefficients in $V_\rho$.
Recently, in \cite[Theorem 5.7]{Su2018}, Su proved that actually we have
\[
H^\bullet({\rm Sh}^\tau,[\mathcal{V}_\rho]^{\rm can}) = H^\bullet(\frak{P}_\tau,K_\infty; \mathcal{A}(\GU_n(\A),{\sf w})\otimes V_\rho),
\]
where $\mathcal{A}(\GU_n(\A),{\sf w})$ is the space of automorphic forms $f$ on $\GU_n(\A)$ such that $f(ag) = a^{{\sf w}}f(g)$ for all $a \in \R_+$ and $g \in \GU_n(\A)$.
Let $\mathcal{A}_0(\GU_n(\A),{\sf w})$ and $\mathcal{A}_{(2)}(\GU_n(\A),{\sf w})$ be the subspaces of $\mathcal{A}(\GU_n(\A),{\sf w})$ consisting of cusp forms and square-integrable automorphic forms respectively.
The natural inclusions then induce $\GU_n(\A_f)$-equivariant homomorphisms
\begin{align*}
f_0^\bullet : H^\bullet(\frak{P}_\tau,K_\infty; \mathcal{A}_0(\GU_n(\A),{\sf w})\otimes V_\rho) \longrightarrow H^\bullet({\rm Sh}^\tau,[\mathcal{V}_\rho]^{\rm can}),\\
f_{(2)}^\bullet : H^\bullet(\frak{P}_\tau,K_\infty; \mathcal{A}_{(2)}(\GU_n(\A),{\sf w})\otimes V_\rho) \longrightarrow H^\bullet({\rm Sh}^\tau,[\mathcal{V}_\rho]^{\rm can}).
\end{align*}
By \cite[Proposition 3.6 and Theorem 5.3]{Harris1990}, we have
\begin{align}\label{E:rational structure pf 1}
{\rm kernel}(f_0^\bullet) = 0,\quad  {\rm image}(f_0^\bullet) \subseteq H_!^\bullet({\rm Sh}^\tau,[\mathcal{V}_\rho]) \subseteq {\rm image}(f_{(2)}^\bullet). 
\end{align}
In particular, $H_!^\bullet({\rm Sh}^\tau,[\mathcal{V}_\rho])$ is a semisimple $\GU_n(\A_f)$-module.

\subsection{Rational structures}\label{SS:rational structure}

Let $\itPi$ be a globally generic cuspidal automorphic representation of $\GU_n(\A)$. We assume that $\itPi_\infty$ is a $C$-algebraic discrete series representation. 
We write $\itPi_\infty^\iota = \itPi_\infty$ and $\itPi_\infty^{\overline{\iota}} = \itPi_\infty^c$. 
For $\sigma \in {\rm Aut}(\C)$, we define the $\sigma$-conjugate of $\itPi$ by
\[
{}^\sigma\!\itPi = \itPi_\infty^{\sigma\circ\iota}\otimes{}^\sigma\!\itPi_f.
\]
We denote by $\Q(\itPi)$ the rationality field of $\itPi$, that is, the fixed field in $\C$ of the set of automorphisms $\sigma \in {\rm Aut}(\C)$ such that ${}^\sigma\!\itPi \cong \itPi$.

In the following proposition, we show that being globally generic is an arithmetic property. The proof relies on the result of Harris on coherent cohomology of Shimura varieties and Arthur's endoscopic classification for $\U_n(\A)$ due to Mok \cite{Mok2015} subject to some results in unpublished preprints of Arthur. Recently, Atobe--Gan--Ichino--Kaletha--M\'{i}nguez--Shin \cite{AGIKMS2024} have address these issues, making the endoscopic classification conditional only on the twisted weighted fundamental lemma.

\begin{prop}\label{P:rational structure}
The following assertions hold.
\begin{itemize}
\item[(1)] For $\sigma \in {\rm Aut}(\C)$, the representation ${}^\sigma\!\itPi$ is cuspidal automorphic and globally generic.
\item[(2)] $\itPi_f$ is defined over $\K\cdot\Q(\itPi_f)$.
\item[(3)] We have $\Q(\itPi) = \Q(\itPi_f)$.
\end{itemize}
\end{prop}

\begin{proof} 
By assumption, $\itPi_\infty = \itPi_\Lambda$ for some admissible $R_{\iota,c}^+$-positive weight $\Lambda = (\lambda;{\sf w})$. We may assume $\Lambda$ is of the form (\ref{E:generic parameter}). Let $V_\rho = V_{\Lambda^\vee-\delta_\iota}$ be the irreducible algebraic representation of $K_\infty$ of $R_{\iota,c}^+$-highest weight $\Lambda^\vee-\delta_\iota$.
We denote by $H_!^\bullet({\rm Sh}^\tau,[\mathcal{V}_\rho])[\itPi_f]$ the $\itPi_f$-isotypic component of the interor cohomology. 
Since the $(\frak{P}_\iota,K_\infty)$-cohomology of $\itPi_\infty \otimes V_\rho$ is non-zero at degree $q(\Lambda)$ according to the discussion in \S\,\ref{SS:d.s. rep.}, by (\ref{E:rational structure pf 1}) the $\itPi_f$-isotypic part is non-zero at degree $\bullet = q(\Lambda)$.
Let $\sigma \in {\rm Aut}(\C)$. Put
\[
{}^\sigma\!\Lambda = \begin{cases}
\Lambda & \mbox{ if $\sigma \vert_\K = {\rm id}$},\\
\Lambda^c & \mbox{ if $\sigma \vert_\K = c$}.
\end{cases}
\]
Thus $\itPi_\infty^{\sigma\circ\iota} \cong \itPi_{{}^\sigma\!\Lambda}$ is the discrete series representation of Harish-Chandra parameter ${}^\sigma\!\Lambda$ in the notion of \S\,\ref{SS:d.s. rep.}.
By the $\GU_n(\A_f)$-equivariance, we have
\begin{align}\label{E:rational structure pf 0}
T_\sigma ( H^\bullet_!({\rm Sh}^\iota,[\mathcal{V}_\rho])[\itPi_f] ) =  H^\bullet_!({\rm Sh}^{\sigma\circ\iota},[\mathcal{V}_{{}^\sigma\!\rho}])[{}^\sigma\!\itPi_f].
\end{align}
In particular, the ${}^\sigma\!\itPi_f$-isotypic component on the right-hand side is non-zero at degree $q(\Lambda)$.
By (\ref{E:rational structure pf 1}), there exists an irreducible discrete automorphic representation $\itPi'$ of $\GU_n(\A)$ such that $\itPi'_f \cong {}^\sigma\!\itPi_f$ and
\begin{align}\label{E:rational structure pf 2}
H^{q(\Lambda)}(\frak{P}_{\sigma\circ\iota},K_\infty; \itPi'_\infty \otimes V_{{}^\sigma\!\rho}) \neq 0.
\end{align}
Note that $\Lambda^c$ is $R_{\overline{\iota},c}^+$-positive, $q(\Lambda^c) = q(\Lambda)$, and $V_{\rho^c} = V_{(\Lambda^c)^\vee-\delta_{\overline{\iota}}}$. Therefore, by (\ref{E:rational structure pf 2}) and \cite[Proposition 4.3.2]{Harris1990}, the infinitesimal characters of $\itPi'_\infty$ and $\itPi_{{}^\sigma\!\Lambda}$ are equal.
%Fix $\varepsilon \in \{\pm1\}$ such that $\pi_{{}^\sigma\!\Lambda}$ is $\psi_{N_n,\infty}^\varepsilon$-generic. Note that the choice of $\varepsilon$ is unique if $n$ is even. It depends on $g_\infty$ in (\ref{E:polarization}) (cf.\,\cite[Proposition A.3]{Atobe2020}).
%By \cite[Lemma 3.5]{LM2015}, $\itPi' \vert_{\U_n(\A)}$ (restriction of automorphic forms to $\U_n(\A)$) admits a unique globally $\psi_{N_n}^\varepsilon$-generic irreducible constitute.
Let $\itPi_1'$ be an irreducible constitute of $\itPi' \vert_{\U_n(\A)}$ (restriction of automorphic forms to $\U_n(\A)$) such that 
\begin{align}\label{E:rational structure pf 3}
H^{q(\Lambda)}(\frak{P}_{\sigma\circ\iota},K_\infty; \itPi'_{1,\infty} \otimes V_{{}^\sigma\!\rho}) \neq 0.
\end{align}
With respect to the Whittaker datum $(B_n \cap \U_n,\psi_{N_n})$, we have Arthur's endoscopic classification of discrete automorphic spectrum of $\U_n(\A)$ established by Mok \cite{Mok2015}.
Consider the functorial lift of $\itPi_1'$ to $\GL_n(\A_\K)$ with respect to the standard representation
${}^L{\rm U}_n \rightarrow {}^L({\rm Res}_{\K/\Q}\GL_n)$. The functorial lift is an isobaric automorphic representation which we denote by $\itPsi = \itPsi(\itPi_1)$.
For each prime number $p$, $\itPi_{1,p}'$ appears in $\itPi'_p\vert_{\U_n(\Q_p)} \cong {}^\sigma\!\itPi_p\vert_{\U_n(\Q_p)}$. 
Thus $\itPi_{1,p}'$ is generic.
In particular, $\itPi_1'$ is almost locally generic. It follows that the cuspidal summands of $\itPsi$ are all unitary. In other words, the global Arthur parameter of $\itPi_1'$ is generic. 
In this case, the associated local Arthur packets are local $L$-packets. 
The condition on the infinitesimal character of $\itPi_{1,\infty}'$ then implies that $\itPsi_\infty$ is also $C$-algebraic, that is, algebraic in the original definition of Clozel in \cite[Definition 1.8]{Clozel1990}. It then follows from the purity lemma \cite[Lemme 4.9]{Clozel1990} that $\itPsi_\infty$ is tempered. Hence $\itPi_{1,\infty}'$ is also tempered. By (\ref{E:rational structure pf 3}), \cite[Theorem 3.5]{Harris1990b}, and the discussion in \S\,\ref{SS:d.s. rep.}, we then conclude that $\itPi_{1,\infty}' \cong \pi_{{}^\sigma\!\Lambda}\vert_{\U_n(\R)}$.
Since $\itPi_{1,\infty}'$ appears in $\itPi_\infty' \vert_{\U_n(\R)}$, the irreducibility of $\itPi_\infty'$ then implies that $\itPi_\infty' \cong \itPi_{{}^\sigma\!\Lambda}$. 
Therefore, ${}^\sigma\!\itPi \cong \itPi'$ is discrete automorphic. Since $\itPi_{{}^\sigma\!\Lambda}$ is essentially tempered, this implies that $\itPi'$ is cuspidal by the result of Wallach \cite[Theorem 4.3]{Wallach1984}. To prove the first assertion, it remains to show that $\itPi'$ is globally generic and appears with multiplicity one in the discrete spectrum. 
By the Arthur's multiplicity formula, since the Arthur parameter of $\itPi_1'$ is generic, there exists a unique 
(counted with multiplicity) discrete automorphic representation $\itPi_1''$ of $\U_n(\A)$ in the global Arthur packet of $\itPi_1'$ such that $\itPi_{1,v}''$ is $\psi_{N_n,v}$-generic for all places $v$. Note that $\itPi_1''$ must be cuspidal by $loc.$ $cit.$.
On the other hand, the automorphic descent of $\itPsi$ to $\U_n(\A)$ is cuspidal, globally $\psi_{N_n}$-generic, and belongs to the Arthur packet of $\itPi_1'$ (cf.\,\cite[Theorem 14]{Soudry2005}).
Therefore, $\itPi_1''$ is the descent of $\itPsi$ and thus globally $\psi_{N_n}$-generic.
By the result of Labesse--Schwermer \cite[Theorem 5.2.2]{LS2019} (with $G=\GU_n$ and $H=\U_n$), $\itPi_1''$ appears in the restriction of some cuspidal automorphic representation $\itPi''$ of $\GU_n(\A)$. Let $v$ be a place of $\Q$. Any $L$-packet of $\U_n(\Q_v)$ contains at most two generic representations (cf.\,\cite[\S\,3]{Kaletha2013}). More precisely, if an $L$-packet contains two generic representations, then the two representations are conjugate to each other by an element of $\GU_n(\Q_p)$ with similitude not in $N_{\K_v/\Q_v}(\K_v^\times)$.
Since $\itPi_{1,v}'$ and $\itPi_{1,v}''$ are in the same $L$-packet, $\itPi_{1,v}'$ appears in $\itPi_v'\vert_{\U_n(\Q_v)}$, and $\itPi_{1,v}''$ appears in $\itPi_v''\vert_{\U_n(\Q_v)}$, we deduce that $\itPi_v' \cong \itPi_v''$. Hence $\itPi' \cong \itPi''$.
By the Arthur's multiplicity formula, the discrete spectrum of $\U_n(\A)$ satisfies multiplicity-one property.
Therefore, $\itPi_1''$ also appears in $\itPi' \vert_{\U_n(\A)}$, which implies that $\itPi'$ is globally generic.
Note that the multiplicity-one property holds for globally generic cuspidal automorphic representations of $\GU_n(\A)$ by \cite[Theorem 5.4.1]{LS2019} (with $G=\GU_n$ and $H=\U_n$). Therefore, actually we have $\itPi' = \itPi''$.

By the reasoning in the above paragraph with $\sigma = {\rm id}$, we see that
\begin{align*}%\label{E:multi. one}
H^\bullet_!({\rm Sh}^\iota,[\mathcal{V}_\rho])[\itPi_f] \cong \itPi_f
\end{align*}
as $\GU_n(\A_f)$-modules. For all $\sigma \in {\rm Aut}(\C/\K\cdot \Q(\itPi_f))$, by (\ref{E:rational structure pf 0}) the isotypic component $H^\bullet_!({\rm Sh}^\iota,[\mathcal{V}_\rho])[\itPi_f]$ is stable under $T_\sigma$. A $\K\cdot \Q(\itPi_f)$-rational structure on the isotypic component is then given by taking the ${\rm Aut}(\C/\K\cdot \Q(\itPi_f))$-invariant subspace (cf.\,\cite[Lemme 3.2.1]{Clozel1990})
\[
\left(H^\bullet_!({\rm Sh}^\iota,[\mathcal{V}_\rho])[\itPi_f]\right)^{{\rm Aut}(\C/\K\cdot \Q(\itPi_f))}.
\]
Therefore, $\itPi_f$ is defined over $\K\cdot \Q(\itPi_f)$. 
Now we prove the third assertion. It is clear that ${\rm Aut}(\C/\Q(\itPi)) \subseteq {\rm Aut}(\C/\Q(\itPi_f))$. Conversely, let $\sigma \in {\rm Aut}(\C/\Q(\itPi_f))$.
Then we have $H^\bullet_!({\rm Sh}^{\sigma\circ\iota},[\mathcal{V}_{{}^\sigma\!\rho}])[\itPi_f] \neq 0$ by (\ref{E:rational structure pf 0}). The arguments in the first paragraph show that $\itPi':=\itPi_\infty^{\sigma\circ\iota}\otimes \itPi_f$ is cuspidal automorphic and globally generic.
Let $\itPi_1$ and $\itPi_1'$ be the unique globally $\psi_{N_n}$-generic irreducible constitute of $\itPi \vert_{\U_n(\A)}$ and $\itPi_1'\vert_{\U_n(\A)}$ respectively. 
Then $\itPi_1$ and $\itPi_1''$ are nearly equivalent. By Arthur's multiplicity formula, we deduce that $\itPi_1 = \itPi_1''$. Thus $\itPi_{1,\infty}$ appears in both $\itPi_\infty^{\iota} \vert_{\U_n(\R)}$ and $\itPi_\infty^{\sigma\circ\iota} \vert_{\U_n(\R)}$.
This implies that $\itPi_\infty^{\iota} \cong \itPi_\infty^{\sigma\circ\iota}$. 
This completes the proof.
\end{proof}

In the following discussion, irreducible algebraic representations of $K_\infty$ are parametrized by $R_{\iota,c}^+$-dominant integral weights.
For $\tau \in \{\iota,\overline{\iota}\}$, we fix a highest weight and a lowest weight Whittaker functions
\[
 W_{\infty,+}^{\tau} \in \mathcal{W}(\itPi_\infty^\tau,\psi_{N_n,\infty}),\quad
 W_{\infty,-}^{\tau} \in \mathcal{W}((\itPi_\infty^\vee)^\tau,\psi_{N_n,\infty}^{-1})
\]
respectively in the minimal $K_\infty$-type such that $W_{\infty,+}^\tau(1) \neq 0$ and $W_{\infty,-}^\tau(1) \neq 0$.
For $\sigma \in {\rm Aut}(\C)$, we denote by $u_\sigma \in \widehat{\Z}^\times$ the unique element such that $\sigma(\psi(x)) = \psi(u_\sigma x)$ for all $x \in \A_f$.
Let $d_\sigma \in \GU_n(\A_f)$ be the diagonal matrix defined so that its $(i,i)$-entry is equal to
\[
\begin{cases}
u_\sigma^{(-2+2i-n)/2} & \mbox{ if $n$ is even and $1 \leq i \leq n/2$},\\
u_\sigma^{n-i} & \mbox{ if $n$ is even and $n/2 < i \leq n$},\\
u_\sigma^{(-1+2i-n)/2} & \mbox{ if $n$ is odd}.
\end{cases}
\]
We then have
\begin{align}\label{E:sigma on psi}
\sigma (\psi_{N_n}(u)) = \psi_{N_n}(d_\sigma^{-1} u d_\sigma),\quad u \in N_n(\A_f),\,\sigma\in{\rm Aut}(\C).
\end{align}
Recall ${}^\sigma\!\itPi_+$ is the space of cusp forms in ${}^\sigma\!\itPi$ whose archimedean components are highest weight vectors in the minimal $K_\infty$-type of ${}^\sigma\!\itPi_\infty$. With respect to the choice $W_{\infty,+}^{\sigma\circ\iota}$, we have a $\GU_n(\A_f)$-equivariant isomorphism between $\mathcal{W}({}^\sigma\!\itPi_f,\psi_{N_n,f})$ and ${}^\sigma\!\itPi_+$ defined as in (\ref{E:Whittaker identification}).
Let 
\[
t_\sigma : \mathcal{W}(\itPi_f,\psi_{N_n,f}) \longrightarrow \mathcal{W}({}^\sigma\!\itPi_f,\psi_{N_n,f})
\]
be the $\sigma$-linear $\GU_n(\A_f)$-equivariant isomorphism defined by 
\[
t_\sigma W(g) = \sigma\left(W( d_\sigma g) \right),\quad g \in \GU_n(\A).
\]
We then have a $\sigma$-linear $\GU_n(\A_f)$-equivariant isomorphism
\begin{align}\label{E:global sigma-linear}
\itPi_+ \longrightarrow {}^\sigma\!\itPi_{+},\quad  f \longmapsto {}^\sigma\!f
\end{align}
defined so that the following diagram is commutative
\[
\begin{tikzcd}[row sep=normal, column sep=normal]
\itPi_+\arrow[r]&{}^\sigma\!\itPi_+\\
\mathcal{W}(\itPi_f,\psi_{N_n,f})\arrow[r, "t_\sigma"]\arrow[u] &\mathcal{W}({}^\sigma\!\itPi_f,\psi_{N_n,f})\arrow[u].
\end{tikzcd}
\]
For a subfield $E$ of $\C$ containing $\Q(\itPi)$, which is equal to $\Q(\itPi_f)$ by Proposition \ref{P:rational structure}-(3), we say $f \in \itPi_+$ is Whittaker-rational over $E$ if and only if ${}^\sigma\!f = f$ for all $\sigma \in {\rm Aut}(\C/E)$.
Similarly, with respect to the choice $W_{\infty,-}^{\sigma\circ\iota}$, we can identify $\mathcal{W}({}^\sigma\!\itPi_f^\vee,\psi_{N_n,f}^{-1})$ with ${}^\sigma\!\itPi_-^\vee$, and we have the notion of been Whittaker-rational for cusp forms in $\itPi_-^\vee$.

For arithmetic application, it is customary to choose normalized newforms as specific Whittaker-rational cusp forms. 
The notion of newform depends on the distinguished open compact subgroup. For instance, the local newform theory is known for $\GL_N$ and $\GSp_4$. Recently, Atobe--Oi--Yasuda \cite{AOY2024} and Cheng \cite{Cheng2023} have carry out the case for unramified odd unitary groups. In particular, normalized newforms (hence Whittaker-rational cusp forms over $\Q(\itPi)$) exist for $\itPi_+$ under the following conditions:
\begin{itemize}
\item If $p \mid D_\K$, then $\itPi_p$ is spherical. 
\item If $n$ is even, then $\itPi_p$ is unramified for all prime $p$ inert in $\K$.
\end{itemize}
By Proposition \ref{P:rational structure}-(2), $\itPi_f$ is defined over $\K\cdot\Q(\itPi_f)$. If $E$ contains $\K\cdot\Q(\itPi_f)$, then the existence of a non-zero Whittaker-rational cusp form over $E$ in $\itPi_+$ actually implies that the space of Whittaker-rational cusp forms over $E$ defines an $E$-rational structure on $\itPi_+$. 
More generally, we have the following lemma.

\begin{lemma}\label{L:rational structure}
Let $(\pi,V_\pi)$ be an irreducible complex representation of a group $G$. Assume the following conditions are satisfied:
\begin{itemize}
\item[(1)] The representation $(\pi,V_\pi)$ is defined over a subfield $E$ of $\C$. 
\item[(2)] For each $\sigma \in {\rm Aut}(\C/E)$, we have a $\sigma$-linear $G$-equivariant isomorphism $t_\sigma : V_\pi \rightarrow V_\pi$. %such that $t_{\sigma_1\circ\sigma_2} = t_{\sigma_1}\circ t_{\sigma_2}$.
\item[(3)] Schur's lemma holds for $\pi$.
\end{itemize} 
Then $V_\pi^{{\rm Aut}(\C/E)} := \{v \in V_\pi \,\vert\, t_\sigma v = v \mbox{ for all }\sigma \in {\rm Aut}(\C/E)\}$ defines an $E$-rational structure of $(\pi,V_\pi)$ if and only if it is non-zero. 
\end{lemma}

\begin{proof}
Let $W \subset V_\pi$ be an $E$-rational structure. For $\sigma \in {\rm Aut}(\C/E)$, define a $\sigma$-linear $G$-equivariant isomorphism $t_\sigma' : V_\pi \rightarrow V_\pi$ by 
\[
t_\sigma' \left(\sum_i z_iv_i \right) = \sum_i\sigma(z_i)v_i,\quad z_i \in \C,\,v_i \in W.
\]
Then $t_\sigma' \circ t_\sigma^{-1}$ is $G$-equivariant linear isomorphism. By Schur's lemma, there exists $\lambda_\sigma \in \C^\times$ such that
\[
t_\sigma' \circ t_\sigma^{-1} = \lambda_\sigma \cdot {\rm id}.
\]
In other words, 
\[
t_\sigma\left(\sum_i z_iv_i \right) = \lambda_\sigma^{-1}\cdot \sum_i\sigma(z_i)v_i,\quad z_i \in \C,\,v_i \in W.
\]
Let $v  \in V_\pi$ be a non-zero vector. Write $v = \sum_{i=1}^r z_iv_i$ for some $z_i \in \C$ and $v_i \in W$. We may assume $v_1,\cdots,v_r$ are linearly independent over $E$ (hence over $\C$) and $z_i \neq 0$ for all $i$. Then $v\in V_\pi^{{\rm Aut}(\C/E)}$ if and only if 
\begin{align}\label{E:rational structure pf}
\lambda_\sigma = \frac{\sigma(z_1)}{z_1} = \cdots = \frac{\sigma(z_r)}{z_r},\quad  \sigma \in {\rm Aut}(\C/E).
\end{align}
If $V_\pi^{{\rm Aut}(\C/E)}$ is non-zero, the above condition then implies that there exists $z \in \C^\times$ such that $\lambda_\sigma = \tfrac{\sigma(z)}{z}$ for all $\sigma \in {\rm Aut}(\C/E)$.
In this case, by using (\ref{E:rational structure pf}) we then obtain $V_\pi^{{\rm Aut}(\C/E)} = z\cdot W$.
This completes the proof.
\end{proof}

\subsection{Adjoint $L$-values for unitary similitude groups}

We keep the notation and assumption as in \S\,\ref{SS:rational structure}.
In the definition of adjoint $L$-functions for a connected reductive group $G$ over $\Q$, ${\rm Ad}$ is the adjoint action of its Langlands dual group ${}^LG$ on ${\rm Lie}(\widehat{G})$. Let ${\rm Ad}_1$ be the adjoint action of ${}^LG$ on ${\rm Lie}(\widehat{G}) / Z({\rm Lie}(\widehat{G}))^{\Gal(\overline{\Q}/\Q)}$.
In the case $G=\GU_n$ considered in this section, we then have
\[
L(s,\itPi,{\rm Ad}) = \zeta(s) L(s,\itPi,{\rm Ad}_1),
\]
where $\zeta(s)$ is the completed Riemann zeta function. 
Note that $L(s,\itPi,{\rm Ad}_1)$ is holomorphic and non-vanishing at $s=1$.
Also we have
\[
\Delta_{\GU_n}(s) = \zeta(s)\prod_{i=1}^nL(s+i-1,\omega_{\K/\Q}^i),
\]
where $\omega_{\K/\Q}$ is the quadratic Hecke character of $\A^\times$ associated to $\K/\Q$ by class field theory.
In the literature, Arthur's conjecture for $\GU_n(\A)$ is not yet proved. Nonetheless, we can deduce from \cite[\S\,Lemma 3.5]{LM2015} that Conjecture \ref{C:LM} still holds for $\itPi$ on the automorphic side, that is, (\ref{E:LM}) holds for $\itPi$ with the proportional constant given by certain automorphic invariant of $\itPi$.

\begin{prop}\label{P: LM for GU}
Let $k_{\itPi}$ be the number of cuspidal summands in the functorial lift of $\itPi_1$ to $\GL_n(\A_\K)$, where $\itPi_1$ is the unique globally $\psi_{N_n}$-generic cuspidal automorphic representation of $\U_n(\A)$ appears in $\itPi \vert_{\U_n(\A)}$ (restriction of cusp forms to $\U_n(\A)$). Let $X(\itPi)$ be the set of Hecke characters $\chi$ of $\A^\times$ such that $\itPi\otimes(\chi\circ\nu) \cong \itPi$. Then (\ref{E:LM}) holds for $\itPi$ with $|\mathcal{S}_\itPi|$ replaced by $2^{k_{\itPi}-1}|X(\itPi)|^{-1}$.
\end{prop}

\begin{proof}
By the results of Beuzart-Plessis--Chaudouard \cite{BPC2023} and Morimoto \cite{Morimoto2024}, the Lapid--Mao conjecture holds for $\U_n(\A)$. In particular, (\ref{E:LM}) holds for $\itPi_1$ with $|\mathcal{S}_{\itPi_1}| = 2^{k_\itPi-1}$. 
On the other hand, as explained in the end of the first paragraph of the proof of Proposition {\ref{P:rational structure}}, the multiplicity of $\itPi$ in the cuspidal spectrum is one. Therefore, it follows from \cite[Lemma 3.5]{LM2015} that (\ref{E:LM}) holds for $\itPi$ with proportional constant $|\mathcal{S}_{\itPi_1}|^{-1}|X(\itPi)|$.
\end{proof}

With respect to the choices $W_{\infty,+}^\tau$ and $W_{\infty,-}^\tau$ for $\tau \in \{\iota,\overline{\iota}\}$, we have a non-zero constant $C_{N_n,\itPi_\infty^{\tau}}$ defined as in (\ref{E:C_infty 1}).
Following result is a refinement of Theorem \ref{T: main 0} for unitary similitude groups.

\begin{thm}\label{T: main 1}
%Assume $\itPi$ appears in the cuspidal spectrum of $\GU_n(\A)$ with multiplicity one when $n$ is even.
Let $\itPi$ be a globally generic cuspidal automorphic representation of $\GU_n(\A)$. Assume that the $\itPi$ is $C$-algebraic and $\itPi_\infty$ is a discrete series representation.  
For $f \in \itPi_+$ and $f' \in \itPi^\vee_-$ with $\<f,f'\>\neq 0$, we have
\begin{align*}
\sigma \left( \frac{L^S(1,\itPi,{\rm Ad}_1)}{\pi^{n(n+1)/2}\cdot \sqrt{D_\K}^{\lfloor (n+1)/2 \rfloor} \cdot C_{N_n,\itPi_\infty^\iota}\cdot \<f,f'\>} \right) = \frac{L^S(1,{}^\sigma\!\itPi,{\rm Ad}_1)}{\pi^{n(n+1)/2}\cdot \sqrt{D_\K}^{\lfloor (n+1)/2 \rfloor} \cdot C_{N_n,\itPi_\infty^{\sigma\circ\iota}}\cdot \<{}^\sigma\!f,{}^\sigma\!f'\>}
\end{align*}
for all $\sigma \in {\rm Aut}(\C)$,
where $S$ is a sufficiently large finite set of places of containing the archimedean place.
In particular, if $f$ and $f'$ are Whittaker-rational over $\Q(\itPi_f)$, then we have
\[
\frac{L^S(1,\itPi,{\rm Ad}_1)}{\pi^{n(n+1)/2}\cdot \sqrt{D_\K}^{\lfloor (n+1)/2 \rfloor} \cdot C_{N_n,\itPi_\infty^\iota}\cdot \<f,f'\>} \in \Q(\itPi_f).
\]
\end{thm}

\begin{proof}
We follow the notation in \S\,\ref{SS: main 0 pf} with $G=\GU_n$.
Fix $\sigma \in {\rm Aut}(\C)$.
For each prime number $p$, we have the following refinement of Lemma \ref{L:Galois equiv. LM 0}, which can be proved in a similar way using (\ref{E:sigma on psi}):
For $f_p \in \itPi_p$ and $f_p' \in \itPi_p^\vee$, we have
\begin{align}\label{E:main 1 pf 1}
\sigma(I_p(f_p\otimes f_p')) = I_{p,\sigma}\left({}^\sigma\!\itPi_p(d_{\sigma,p}^{-1}){}^\sigma\!f_p\otimes{}^\sigma\!\itPi_p^\vee(d_{\sigma,p}^{-1}){}^\sigma\!f_p'\right).
\end{align}
By Proposition \ref{P: LM for GU}, (\ref{E:main 1 pf 1}), and proceed as in the proof of Theorem \ref{T: main 0}, we then have
\[
\sigma \left(\left.\left(\frac{L^S(s,\itPi,{\rm Ad})}{\Delta_{\GU_n}^S(s)}\right)\right\vert_{s=1}\frac{ 2^{k_{\itPi}-1}|X(\itPi)|^{-1}}{C_{N_n,\itPi_\infty^\iota}\cdot\<f,f'\>} \right) = \left.\left(\frac{L^S(s,{}^\sigma\!\itPi,{\rm Ad})}{\Delta_{\GU_n}^S(s)}\right)\right\vert_{s=1}\frac{2^{k_{{}^\sigma\!\itPi}-1}|X({}^\sigma\!\itPi)|^{-1}}{C_{N_n,\itPi_\infty^{\sigma\circ\iota}}\cdot\<{}^\sigma\!f,{}^\sigma\!f'\>}
\]
for any $f \in \itPi_+$ and $f' \in \itPi^\vee_-$ with $\<f,f'\>\neq 0$.
It is clear that $\chi \in X(\itPi)$ if and only if ${}^\sigma\!\chi \in X({}^\sigma\!\itPi)$ (in fact, $\chi= {\bf 1}$ or $\omega_{\K/\Q}$). In particular, $|X(\itPi)| = |X({}^\sigma\!\itPi)|$. 
Now we show that $k_\itPi = k_{{}^\sigma\!\itPi}$.
Let $\itPi_1$ and $({}^\sigma\!\itPi)_1$ be the unique globally $\psi_{N_n}$-generic cuspidal constitute of $\itPi\vert_{\U_n(\A)}$ and $({}^\sigma\!\itPi)\vert_{\U_n(\A)}$ respectively. Consider their functorial lifts $\itPsi(\itPi_1)$ and $\itPsi(({}^\sigma\!\itPi)_1)$ to $\GL_n(\A_\K)$, which are isobaric automorphic representations.
We have
\[
\itPsi(\itPi_1) = \itSigma_1 \boxplus \cdots \boxplus \itSigma_k
\]
for some conjugate self-dual cuspidal automorphic representation $\itSigma_i$ of $\GL_{n_i}(\A_\K)$ for $1 \leq i \leq k = k_\itPi$.
Since $\itPi_{1,\infty}$ is a discrete series representation of $\U_n(\R)$, its functorial lift $\itPsi(\itPi_1)_\infty$ is regular algebraic (cf.\,\cite[Theorem 5.6-(4)]{RS2018}). This implies that $\itSigma_i \otimes |\mbox{ }|_\A^{(n-n_i)/2}$ is regular algebraic for all $i$. It then follows from the result of Clozel \cite[Th\'eor\`eme 3.19]{Clozel1990} that ${}^\sigma\!(\itSigma_i \otimes |\mbox{ }|_\A^{(n-n_i)/2})$ is cuspidal automorphic for all $i$ and we have
\[
{}^\sigma\!\itPsi(\itPi_1) = \bigboxplus_{i=1}^k{}^\sigma\!(\itSigma_i \otimes |\mbox{ }|_\A^{(n-n_i)/2})\otimes |\mbox{ }|_\A^{-(n-n_i)/2}.
\]
Therefore, if ${}^\sigma\!\itPsi(\itPi_1) = \itPsi(({}^\sigma\!\itPi)_1)$, then $k_\itPi = k_{{}^\sigma\!\itPi}$.
For each prime number $p$, it is clear that ${}^\sigma\!\itPi_{1,p}$ and $({}^\sigma\!\itPi)_{1,p}$, which are the local components at $p$ of ${}^\sigma\!\itPi_1$ and $({}^\sigma\!\itPi)_1$ respectively, both appear in ${}^\sigma\!\itPi_p \vert_{\U_n(\Q_p)}$. Hence ${}^\sigma\!\itPi_{1,p}$ is isomorphic to $({}^\sigma\!\itPi)_{1,p}\circ {\rm Ad}(g_p)$ for some $g_p \in \GU_n(\Q_p)$. In particular, if $p \notin S$, then one can easily verify that ${}^\sigma\!\itPi_{1,p}$ and $({}^\sigma\!\itPi)_{1,p}$ have the same functorial lift to $\GL_n(\Q_p)$ as their Satake parameters are the same. On the other hand, by \cite[Lemma 9.2]{GR2013}, the functorial lift of ${}^\sigma\!\itPi_{1,p}$ is isomorphic to ${}^\sigma\!\itPsi(\itPi_1)_p$ for all $p \notin S$. We thus conclude that ${}^\sigma\!\itPsi(\itPi_1)$ and $\itPsi(({}^\sigma\!\itPi)_1)$ are nearly equivalent. Therefore, ${}^\sigma\!\itPsi(\itPi_1) = \itPsi(({}^\sigma\!\itPi)_1)$ by the strong multiplicity one theorem of Jacquet and Shalika \cite[Theorem 4.4]{JS1981b}.
Finally, note that
\[
\left.\left(\frac{L^S(s,\itPi,{\rm Ad})}{\Delta_{\GU_n}^S(s)}\right)\right\vert_{s=1} = \frac{L^S(1,\itPi,{\rm Ad}_1)}{\prod_{i=1}^nL^S(i,\omega_{\K/\Q}^i)}.
\]
By the algebraicity of critical values of Dirichlet $L$-functions, we have
\[
\prod_{i=1}^nL^S(i,\omega_{\K/\Q}^i) \in \pi^{n(n+1)/2}\cdot \sqrt{D_\K}^{\lfloor (n+1)/2 \rfloor}\cdot \Q^\times.
\]
This completes the proof.
\end{proof}

\section{An explicit formula for ${\rm U}(2,1)$}

\subsection{Basic setting}\label{E:basic}
Let $\K$ be an imaginary quadratic field of discriminant $-D_\K$. Let $\delta = {\sqrt{-D_\K}}$.
We write $x\mapsto x^c$ for the non-trivial automorphism of $\K$.
Let $G$ be a reductive group over $\Q$ defined by
\[
G = \left\{ g \in {\rm Res}_{\K/\Q}\GL_3 \,\left\vert\, {}^t\!g^c Q g = Q  \right.\right\},
\]
where $Q$ is the $c$-hermitian symmetric matrix given by
\[
Q = \bp &&\delta^{-1} \\ &1& \\ -\delta^{-1}&& \ep.
\]
Note that $G$ admits a unique Whittaker datum up to conjugation.
Let $N$ be the maximal unipotent subgroup of $G$ given by
\begin{align}\label{E:additive character}
N = \left.\left\{ n(x,y)=  \bp 1 & -\delta x & -\delta y+\tfrac{1}{2}\delta xx^c \\  & 1 & -x^c \\  &  & 1 \ep\,\right\vert\, x,y \in {\rm Res}_{\K/\Q}\mathbb{G}_a,\,y+y^c=0 \right\}.
\end{align}
Let $\psi_N$ be the non-degenerated character of $N(\Q)\backslash N(\A)$ defined by
\[
\psi_N\left ( n(x,y) \right) = \psi\circ{\rm Tr}_{\K/\Q}({\delta^{-1}x}),
\]
where $\psi$ is the non-trivial additive character of $\Q\backslash\A$ such that $\psi_\infty(x)=e^{2\pi\sqrt{-1}\,x}$.

For a prime $p$ non-split in $\K$, let $K_p$ be the maximal open compact subgroup defined by 
\[
K_p = G(\Q_p)\cap \GL_3(\mathcal{O}_{p}),
\]
where $\mathcal{O}_p$ is the ring of integers of $\K_p$.
Let $p$ be a prime splits in $\K$. With respect to a fixed $\delta_p \in \Z_p^\times$ such that $\delta_p^2 = -D_\K$, we have an isomorphism 
\[
\K_p \longrightarrow \Q_p \oplus \Q_p, \quad x_1+x_2\delta \longmapsto (x_1+x_2\delta_p,x_1-x_2\delta_p).
\]
Under this isomorphism, we have
\[
G(\Q_p) \cong \left\{(g,Q{}^t\!g^{-1}Q^{-1})\,\right\vert\,g \in \GL_3(\Q_p)\}.
\]
Hereafter we identify $G(\Q_p)$ with $\GL_3(\Q_p)$ via the projection $(g_1,g_2) \mapsto g_1$. For $n \geq 0$, let $K_0(p^n\Z_p)$ be an open compact subgroup of $\GL_3(\Q_p)$ defined by
\[
K_0(p^n\Z_p) = \bp \Z_p&\Z_p&\Z_p \\ \Z_p&\Z_p&\Z_p \\ p^n\Z_p&p^n\Z_p&\Z_p \ep \cap \GL_3(\Z_p).
\]
We also write $K_0(\Z_p) = K_p$.
For a positive integer $M$ whose prime divisors are split in $K$, let $K_0(M)$ be the open compact subgroup of $G(\A_f)$ defined by
\[
K_0(M) = \prod_{p \mid M }K_0(M\Z_p)\times\prod_{p\nmid M}K_p.
\]

We identify $G(\R)$ with the real unitary group $\U(2,1)$ in \S\,\ref{SS:d.s. rep.} as follows.
Let 
\[
g_\infty = \bp \tfrac{-1}{\sqrt{2}\delta}&0&\tfrac{1}{\sqrt{2}} \\ 0& 1 & 0 \\ \tfrac{1}{\sqrt{2}\delta}&0&\tfrac{1}{\sqrt{2}} \ep \in \GL_3(\C).
\]
Note that
$
Q = {}^t\overline{g}_\infty {\rm diag}(1,1,-1) g_\infty$.
Thus we have an isomorphism 
\begin{align}\label{E:iso. 1}
G(\R) \longrightarrow \U(2,1),\quad  g \longmapsto g_\infty g g_\infty^{-1}.
\end{align}
Let $K_\infty$ be the maximal compact subgroup of $G(\R)$ such that $g_\infty K_\infty g_\infty^{-1}$ consisting of ${\rm diag}(k_1,k_2)$ with $k_1 \in \U(2)$ and $k_1 \in \U(1)$.
Under the identification (cf.\,\S\,\ref{SS:d.s. rep.})
\[
 {\rm Lie}(\U(1)^3)_\C^* \ni z_1e_1+z_2e_2+z_3e_3 \longmapsto (z_1,z_2,z_3) \in \C^3,
\]
for each $\lambda \in \Z^3$ with $\lambda_1 \neq \lambda_2$ we let $\pi_\lambda$ be the discrete series representation of $G(\R)$ with Harish-Chandra parameter $\lambda$. Note that $\pm(e_1-e_2)$ are the compact roots.
We parametrize the irreducible algebraic representations of $K_\infty$ by the integral weights which are dominant with respect to the compact root $e_1-e_2$.

\subsection{An explicit formula}

Let $\itPi$ be a globally generic cuspidal automorphic representation of $G(\A)$ satisfying the following conditions:
\begin{itemize}
\item If $p$ is non-split in $\K$, then $\itPi_p$ admits non-zero $K_p$-invariant vectors.
\item $\itPi_\infty = \pi_\lambda$ is a discrete series with Harish-Chandra parameter $\lambda = (\alpha_1,\alpha_3;\alpha_2)$ with $\alpha_1 > \alpha_2 > \alpha_3$ (cf.\,(\ref{E:generic parameter})).
\end{itemize}
By the local nowform theory established by Jacquet--Piatetski-Shapiro--Shalika \cite{JPSS1981} and Reeder \cite{Reeder1991} applied to $\GL_3$, there exists a unique positive integer $M_\itPi$ whose prime divisors are split in $K$ such that the $K_0(M_\itPi)$-invariant subspace of $\itPi_f$ is one-dimensional. We call $M_\itPi$ the conductor of $\itPi$.
%Similarly, $\itPi^\vee$ satisfies the above conditions with $\itPi_\infty^\vee = \pi_{\lambda^\vee}$ where $\lambda^\vee = (-\alpha_3,-\alpha_1;-\alpha_2)$. Note that $M_\itPi = M_{\itPi^\vee}$.
%A cusp form $f$ in $\itPi$ (resp.\,$\itPi^\vee$) is called a newform if the following conditions are satisfied:
%\begin{itemize}
%\item $f$ is right invariant by $K_0(M_\itPi)$.
%\item The archimedean component of $f$ is a highest (resp.\,lowest) weight vector in the minimal $K_\infty$-type of $\itPi_\infty$ (resp.\,$\itPi_\infty^\vee$).
%\end{itemize}
%The subspaces of newforms of $\itPi$ or $\itPi^\vee$ are one-dimensional. 
A cusp form $f$ in $\itPi$ is called a newform if the following conditions are satisfied:
\begin{itemize}
\item $f$ is right invariant by $K_0(M_\itPi)$.
\item The archimedean component of $f$ is a highest weight vector in the minimal $K_\infty$-type of $\itPi_\infty$.
\end{itemize}
The subspace of newforms of $\itPi$ is one-dimensional. 
%We denote by $f_\itPi$ and $f_{\itPi^\vee}$ the normalized newforms of $\itPi$ and $\itPi^\vee$ respectively such that
%\begin{align}\label{E:normalized}
%\color{red}W_{\psi_N}(1;f_\itPi) = K_{\alpha_2-\alpha_3}(),\quad W_{\psi_N^{-1}}(1;f_{\itPi^\vee}) = K_{\alpha_2-\alpha_3}().
%\end{align}
Let $f_\itPi$ be the normalized newform of $\itPi$ such that
\begin{align}\label{E:normalized}
W_{\psi_N}({\rm diag}(\delta,1,\delta^{-c})_\infty;f_\itPi) = K_{\alpha_2-\alpha_3}(4\sqrt{2}\pi).
\end{align}
Here $K_\nu(z)$ is the modified Bessel function defined by
\[
K_\nu(z) = \frac{1}{2} \int_0^\infty e^{-z(t+t^{-1})/2}t^{\nu-1}\,dt,\quad {\rm Re}(z)>0.
\]
In the following theorem, we give an explicit formula between the adjoint $L$-value $L(1,\itPi,{\rm Ad})$ and the Petersson norm %$\<f_\itPi,f_{\itPi^\vee}\>$ defined in (\ref{E:Petersson norm}).
%\[
%\<f_\itPi,f_{\itPi^\vee}\> = \int_{G(\Q)\backslash G(\A)}f_\itPi(g)f_{\itPi^\vee}(g)\,dg^{\rm Tam},
%\]
%where $dg^{\rm Tam}$ is the Tamagawa measure. 
\[
\<f_\itPi,f_{\itPi}\> = \int_{G(\Q)\backslash G(\A)}|f_\itPi(g)|^2\,dg^{\rm Tam},
\]
where $dg^{\rm Tam}$ is the Tamagawa measure.

\begin{thm}\label{T:explicit}
We have
\[
\frac{L(1,\itPi,{\rm Ad})}{\Delta_G(1)} = \<f_\itPi,f_{\itPi}\> \cdot D_\K^2\cdot 2^{-k_\itPi + \alpha_1-2\alpha_2+\alpha_3-1}(\alpha_1-\alpha_3+1)\prod_{p \mid M_\itPi} \frac{L(1,\itPi_p,{\rm Ad})}{\zeta_p(3)\cdot L(1,\itPi_{p,{\rm ur}},{\rm Ad})}.
\]
Here $k_\itPi$ is the number of cuspidal summands in the functorial lift of $\itPi$ to $\GL_3(\A_\K)$%, $\omega(D_\K)$ is the number of prime divisors of $D_\K$
, and $\itPi_{p,{\rm ur}}$ is the unramified component of the first nonzero spherical Bernstein--Zelevinsky derivative of $\itPi_p$ for each $p \mid M_\itPi$.
\end{thm}

\begin{proof}
Write $f_\itPi = \otimes_v f_{\itPi_v} \in \itPi = \otimes_v \itPi_v$. %and $f_{\itPi^\vee} = \otimes_v f_v' \in \itPi^\vee = \otimes_v \itPi_v^\vee$.
For each place $v$ of $\Q$, let $\<\cdot,\cdot\>_v$ be a $G(\Q_v)$-equivariant hermitian pairing on $\itPi_v$. 
Let $I_v$ be the local Whittaker functional defined as in (\ref{E:local Whittaker 1}) with respect to $\<\cdot,\cdot\>_v$ and $\psi_{N,v}$.
By the Lapid--Mao conjecture for $G$ proved by Beuzart-Plessis--Chaudouard \cite{BPC2023} and Morimoto \cite{Morimoto2024}, we have
\[
\frac{L(1,\itPi,{\rm Ad})}{\Delta_G(1)} = 2^{-k_\itPi}\cdot \<f_\itPi,f_{\itPi^\vee}\>\cdot |W_{\psi_N}({\rm diag}(\delta,1,\delta^{-1})_\infty;f_\itPi)|^{-2}\cdot\prod_{v \mid \infty M_\itPi D_\K}C_{\itPi_v,\psi_{N,v}},
\]
where
\begin{align}\label{E:local constants}
C_{\itPi_v,\psi_{N,v}} = \frac{L(1,\itPi_v,{\rm Ad})}{\Delta_{G,v}(1)}\cdot \begin{cases}\displaystyle
\frac{I_p(f_{\itPi_p}\otimes f_{\itPi_p})}{\<f_{\itPi_p}, f_{\itPi_p}\>_p} & \mbox{ if $v=p$},\\
\displaystyle
{\frac{I_\infty(\itPi_\infty({\rm diag}(\delta,1,\delta^{-1})_\infty)f_{\itPi_\infty}\otimes \itPi_\infty({\rm diag}(\delta,1,\delta^{-1})_\infty)f_{\itPi_\infty})}{\<f_{\itPi_\infty}, f_{\itPi_\infty}\>_\infty}} & \mbox{ if $v=\infty$}.
\end{cases} 
\end{align}
These local constants are computed in Propositions \ref{P:split case}, \ref{P:non-split case}, and \ref{P:archimedean} below, from which we have
\[
\prod_{v \mid \infty M_\itPi D_\K}C_{\itPi_v,\psi_{N,v}} = D_\K^2\cdot 2^{\alpha_1-2\alpha_2+\alpha_3-1}(\alpha_1-\alpha_3+1)K_{\alpha_2-\alpha_3}(4\sqrt{2}\pi)^2\prod_{p \mid M_\itPi} \frac{L(1,\itPi_p,{\rm Ad})}{\zeta_p(3)\cdot L(1,\itPi_{p,{\rm ur}},{\rm Ad})}.
\]
This completes the proof.
%\[
%C_{\itPi_p} = \frac{\Delta_{G,p}(1)}{L(1,\itPi_p,{\rm Ad})}\cdot\frac{\<f_{\itPi_p}, f_{\itPi_p^\vee}\>_p}{I_p(f_{\itPi_p}\otimes f_{\itPi_p^\vee})}.
%\]
\end{proof}

An immediate consequence of the explicit formula is the positivity of adjoint $L$-value.
\begin{corollary}
The adjoint $L$-value $L^S(1,\itPi,{\rm Ad})$ is positive, where $S$ is the set of prime divisors of $M_\itPi$.
\end{corollary}

\subsection{Non-archimedean computations}

Let $p$ be a prime divisor of $M_\itPi D_\K$.
In this section, we compute the local constant $C_{\itPi_p,\psi_{N,p}}$.
By our assumption, $p \mid M_\itPi$ (resp.\,$p \mid D_\K$) if and only if $p$ is split (resp.\,non-split) in $\K$. 

\subsubsection{The split case}

Assume $p \mid M_\itPi$. Then $\itPi_p$ is a ramified unitary irreducible admissible representation of $\GL_3(\Q_p)$. Under the isomorphism $G(\Q_p) \cong \GL_3(\Q_p)$, $N(\Q_p)$ is identified with the upper triangular unipotent matrices in $\GL_3(\Q_p)$, and we have
\[
\psi_{N,p}(n) = \psi_p(-\delta_p^{-2}n_{12} - \delta_p^{-1}n_{23}).
\]
In particular, $\psi_{N,p}$ is unramified since $\delta_p \in \Z_p^\times$. 
Let $\mathcal{P}_3$ be the mirabolic subgroup of $\GL_3$ consisting of matrices whose last row is equal to $(0,0,1)$.
Let $\<\cdot,\cdot\>_p$ be the $\GL_3(\Q_p)$-equivariant hermitian pairing on $\mathcal{W}(\itPi_p,\psi_{N,p})$ defined by
\[
\<W,W'\>_p:= \int_{N(\Q_p)\backslash \mathcal{P}_3(\Q_p)} W(g)\overline{W'(g)}\,\overline{dg},
\]
where $\overline{dg} = dn \backslash dg$ is the quotient measure on $N(\Q_p)\backslash \mathcal{P}_3(\Q_p)$ with 
\[
{\rm vol}(N(\Q_p)\cap \GL_3(\Z_p),dn) = {\rm vol}(\mathcal{P}_3(\Z_p),dg)=1.
\]
Let $I_p$ be the Whittaker functional on $\itPi_p \otimes \overline{\itPi_p}$ defined as in (\ref{E:local Whittaker 1}) with respect to $\<\cdot,\cdot\>_p$ and $\psi_{N,p}$, where the Haar measures on $N(\Q_p)$ is normalized so that $N(\Q_p)\cap \GL_3(\Z_p)$ has volume $1$.
By the result of Lapid and Mao \cite[Lemma 4.4]{LM2015}, we have
\begin{align}\label{E:LM splitting lemma}
I_p(W\otimes {W'}) =\zeta_p(1)\zeta_p(2)\cdot W(1)\overline{W'(1)},\quad W,W' \in \mathcal{W}(\itPi_p,\psi_{N,p}).
\end{align}

\begin{prop}\label{P:split case}
We have
\[
C_{\itPi_p,\psi_{N,p}} = \frac{L(1,\itPi_p,{\rm Ad})}{\zeta_p(3)\cdot L(1,\itPi_{p,{\rm ur}},{\rm Ad})}.
\]
Here $\itPi_{p,{\rm ur}}$ is the unramified component of the first nonzero spherical Bernstein--Zelevinsky derivative of $\itPi_p$.
\end{prop}

\begin{proof}
Let $W_{\itPi_p} \in \mathcal{W}(\itPi_p,\psi_{N,p})^{K_0(M_\itPi\Z_p)}$ be the Whittaker newform normalized so that $W_{\itPi_p}(1)=1$.
Note that $\Delta_{G,p}(s) = \zeta_p(s)\zeta_p(s+1)\zeta_p(s+2)$. 
It then follows from (\ref{E:LM splitting lemma}) that
\[
C_{\itPi_p,\psi_{N,p}} = \frac{L(1,\itPi_p,{\rm Ad})}{\zeta_p(3)\cdot\<W_{\itPi_p},W_{\itPi_p}\>_p}.
\]
By the explicit formula for $W_{\itPi_p}$ due to Matringe \cite[Corollary 3.2]{Matringe2013}, one can easily verify that
\[
\<W_{\itPi_p},W_{\itPi_p}\>_p = L(1,\itPi_{p,{\rm ur}},{\rm Ad}).
\]
For instance, please refer to \cite[Proposition 6.4]{HY2023} for the computation. 
\end{proof}

\subsubsection{The non-split case}

Assume $p \mid D_\K$. Then $\itPi_p$ admits a non-zero $K_p$-invariant vector. There exists an unramified character $\chi$ of $\K_p^\times$ such that $\itPi_p$ is isomorphic to an irreducible subrepresentation of the induced representation 
\[
{\rm Ind}_{B(\Q_p)}^{G(\Q_p)}(\chi),
\]
where $B$ is the Borel subgroup of $G$ consisting of upper triangular matrices and $\chi$ is regarded as a character of the diagonal maximal torus by $\chi({\rm diag}(a,b,a^{-c})):=\chi(a)$. More precisely, $\itPi_p$ is isomorphic to the subrepresentation generated by the $K_p$-invariant section $f^\circ \in {\rm Ind}_{B(\Q_p)}^{G(\Q_p)}(\chi)$ normalized so that $f^\circ(1) = 1$. 
Note that the modulus character of $B(\Q_p)$ is given by 
\[
{\rm diag}(a,b,a^{-c}) \longmapsto |a|_{\K_p}^2.
\]
Let $dn$ be the Haar measure on $N(\Q_p)$ with ${\rm vol}(N(\Q_p)\cap \GL_3(\mathcal{O}_p),dn)=1$, and 
\[
w =\bp & & -1 \\ &1& \\ 1&&\ep \in K_p.
\]
Let $\<\cdot,\cdot\>_p$ be the $G(\Q_p)$-equivariant hermitian pairing on ${\rm Ind}_{B(\Q_p)}^{G(\Q_p)}(\chi)$ defined as follows. 
If $\chi$ is unitary, then define
\[
\<f,f'\>_p := \int_{N(\Q_p)}f(wn) \overline{f'(wn)}\,dn.
\]
If $\chi = |\mbox{ }|_{\K_p}^s$ for some $-1/2 < s < 1/2$ (complementary series), then define
\[
\<f,f'\>_p := \int_{N(\Q_p)}f(wn) \overline{Mf'(wn)}\,dn,
\]
where $M: {\rm Ind}_{B(\Q_p)}^{G(\Q_p)}(\chi) \rightarrow {\rm Ind}_{B(\Q_p)}^{G(\Q_p)}(\chi^{-1})$ is a fixed non-zero intertwining operator. 
Let $I_p$ be the Whittaker functional on $\itPi_p \otimes \overline{\itPi_p}$ defined as in (\ref{E:local Whittaker 1}) with respect to $\<\cdot,\cdot\>_p$ and $\psi_{N,p}$, where the Haar measure on $N(\Q_p)$ is defined as above.

\begin{rmk}
In fact, $\chi$ must be unitary since $\itPi_\infty$ is a discrete series representation and thus the generalized Ramanujan conjecture holds for $\itPi$ by the result of Shin \cite[Corollary 1.3]{Shin2011}.
Nonetheless, the computation in this section applies to any irreducible generic unitary spherical representation of $G(\Q_p)$.
\end{rmk}

%Fix a uniformizer $\varpi_p$ of $\mathcal{O}_p$. For instance, one can take $\varpi_p = 2^{-1}\delta$.
Let $\frak{P}$ be the maximal ideal of $\mathcal{O}_p$ and 
\[
A = {\rm diag}(\frak{P}\smallsetminus \frak{P}^2,\mathcal{O}_p^\times,\frak{P}^{-1}\smallsetminus \mathcal{O}_p) \cap G(\Q_p).
\]
%For $x \in \K_p$ and $y \in \delta^{-1}\Q_p$, let $z(x,y)=-\delta y + \frac{1}{2}\delta x x^c$ be the upper rightmost entry of $n(x,y)$.
For $m \in \Z$, let $N_m$ be an open compact subgroups of $N(\Q_p)$ defined by
\[
N_m = \left\{ n(x,y) \in N(\Q_p)\,\left\vert\,x \in \frak{P}^{\lfloor \frac{m}{2} \rfloor},\, \delta y \in p^{\lfloor \frac{m+1}{2} \rfloor}\Z_p\right\}\right..
\]
Note that $n(x,y) \in N_m$ if and only if
\begin{align}\label{E:C-S}
x \in \frak{P}^{\lfloor \frac{m}{2} \rfloor},\quad -\delta y + \frac{1}{2}\delta x x^c \in \frak{P}^m.
\end{align}
Indeed, if $x \in \frak{P}^{\lfloor \frac{m}{2} \rfloor}$, then $\frac{1}{2}\delta xx^c \in \frak{P}^m$ since $\frac{1}{2}\delta$ is a uniformizer of $\mathcal{O}_p$. Therefore, $-\delta y + \frac{1}{2}\delta x x^c \in \frak{P}^m$ if and only if $\delta y \in \frak{P}^m \cap \Q_p = p^{\lfloor \frac{m+1}{2} \rfloor}\Z_p$.
In particular, $N_0 = N(\Q_p)\cap\GL_3(\mathcal{O}_p)$.

\begin{lemma}\label{L:C-S}
For $m \in \Z$, we have
\[
w(N_{-m} \smallsetminus N_{-m+1})w^{-1} \subset N_{m}A^{m} w (N_{m} \smallsetminus N_{m+1}).
\]
\end{lemma}

\begin{proof}
For $x \in \K_p$ and $y \in \delta^{-1}\Q_p$, we write $z(x,y)=-\delta y + \frac{1}{2}\delta x x^c$ for the upper rightmost entry of $n(x,y)$.
First we show that $n(x,y) \in N_m \smallsetminus N_{m+1}$ if and only if
\begin{align}\label{E:C-S pf 1}
x \in \frak{P}^{\lfloor \frac{m}{2} \rfloor},\quad z(x,y) \in \frak{P}^{m} \smallsetminus \frak{P}^{m+1}.
\end{align}
Let $n(x,y) \in N_m \smallsetminus N_{m+1}$.
If $x \in \frak{P}^{\lfloor \frac{m+1}{2} \rfloor}$, then we have $z(x,y) \in \frak{P}^m \smallsetminus \frak{P}^{m+1}$ by (\ref{E:C-S}).
If $x \in \frak{P}^{\lfloor \frac{m}{2} \rfloor} \smallsetminus \frak{P}^{\lfloor \frac{m+1}{2} \rfloor}$, then $m$ must be odd. Thus $\frac{1}{2}\delta xx^c \in \frak{P}^m \smallsetminus \frak{P}^{m+1}$ and $\delta y \in p^{(m+1)/2}\Z_p\subset \frak{P}^{m+1}$, which imply that $z(x,y) \in \frak{P}^m \smallsetminus \frak{P}^{m+1}$.
Conversely, (\ref{E:C-S pf 1}) implies $n(x,y) \in N_m \smallsetminus N_{m+1}$ by (\ref{E:C-S}).
The assertion of the lemma then follows from (\ref{E:C-S pf 1}) and the following identity which holds for any $n(x,y) \in N(\Q_p)$ with $z=z(x,y) \neq 0$:
\[
w n(x,y) w^{-1} = \bp 1 & -\delta x z^{-c} & -z^{-1}\\ &1&-x^cz^{-1} \\ &&1 \ep {\rm diag}(-z^{-c},z^{-1}z^c,-z) w \bp 1 & \delta x z^{-1} & -z^{-1}\\ &1&x^cz^{-c} \\ &&1 \ep.
\]
\end{proof}

\begin{prop}\label{P:non-split case}
We have $C_{\itPi_p,\psi_{N,p}} = 1$.
\end{prop}

\begin{proof}
We assume $\chi$ is unitary. The case of complementary series can be proved in a similar way.
Note that ${\rm vol}(N_m,dn) = p^{-{\lfloor \frac{m}{2} \rfloor}-{\lfloor \frac{m+1}{2} \rfloor}} = p^{-m}$ for all $m \in \Z$.
Since 
\[
N(\Q_p) = N_0 \sqcup \bigsqcup_{m=1}^\infty (N_{-m} \smallsetminus N_{-m+1}), 
\]
by Lemma \ref{L:C-S} we have
\[
\<f^\circ,f^\circ\> = 1+ \sum_{m=1}^\infty p^{-2m}{\rm vol}(N_{-m} \smallsetminus N_{-m+1},dn) = 1+p^{-1}.
\]
By \cite[Proposition 2.8]{LM2015}, we have
\[
I_p(f\otimes f') = J^{\psi_{N,p}}(f)\overline{J^{\psi_{N,p}}(f')},
\]
where $J^{\psi_{N,p}}$ is the Jacquet integral defined by
\[
J^{\psi_{N,p}}(f) = \int_{N(\Q_p)}f(wn)\overline{\psi_{N,p}(n)}\,dn.
\]
Note that $\psi_p\circ{\rm Tr}_{\K_p/\Q_p}$ is trivial on $\delta^{-1}\mathcal{O}_p$ but non-trivial on $\delta^{-1}\frak{P}^{-m}$ for all $m \geq 1$.
Therefore, $\psi_{N,p}$ is trivial on $N_0$ but non-trivial on $N_{-m}$ for all $m \geq 1$.
Let $\varpi$ be a uniformizer of $\mathcal{O}_p$.
By Lemma \ref{L:C-S}, we then have
\begin{align*}
J^{\psi_{N,p}}(f^\circ) &= 1+ \sum_{m=1}^\infty \chi(\varpi)p^{-m} \int_{N_{-m}\smallsetminus N_{-m+1}}\overline{\psi_{N,p}(n)}\,dn\\
& = 1-\chi(\varpi)p^{-1}.
\end{align*}
Finally, note that
\begin{align*}
L(s,\itPi_p,{\rm Ad}) = \zeta_p(s)L(s,\chi)L(s,\chi^{-1}),\quad \Delta_{G,p}(s) = \zeta_p(s+1).
\end{align*}
We conclude that
\begin{align*}
C_{\itPi_p,\psi_{N,p}} &= \frac{(1-p^{-1})^{-1}(1-\chi(\varpi)p^{-1})^{-1}(1-\chi^{-1}(\varpi)p^{-1})^{-1}}{(1-p^{-2})^{-1}}\cdot \frac{(1-\chi(\varpi)p^{-1})\overline{(1-\chi(\varpi)p^{-1})}}{1+p^{-1}}\\
& = 1.
\end{align*}
This completes the proof.
\end{proof}

\subsection{Archimedean computations}

The main result of this section is Proposition \ref{P:archimedean} on the explicit computation of $C_{\itPi_\infty,\psi_{N,\infty}}$. 
%For brevity, we drop the subscript $\infty$ in this section.
The idea is to realize $\itPi_\infty$ as a local theta lift from $\U_2(\R)$ to $\U_3(\R)$. A key ingredient is \cite[Proposition 4.1]{Morimoto2024} which relates the local Whittaker functional for $\U_{2n+1}(\R)$ with that of $\U_{2n}(\R)$ under local theta correspondence. The computation of $C_{\itPi_\infty,\psi_{N,\infty}}$ then boils down to the computation of the ratio appears in $loc.$ $cit.$ for discrete series representations in the case $n=1$. For higher $n$, the computation would be an interesting and challenging problem to consider. 

\subsubsection{Weil representation}\label{SS:Weil rep}

Let $(W,(\cdot,\cdot)_W)$ be the skew-hermitian space over $\C$ defined by $W = {\rm M}_{1,2}(\C)$ (row vectors) and 
\[
(w_1,w_2)_W = w_1\bp  & 1 \\ -1 &  \ep {}^t\overline{w_2}.
\]
Let $(V,(\cdot,\cdot)_V)$ be the hermitian space over $\C$ defined by $V = {\rm M}_{3,1}(\C)$ (column vectors) and 
\[
(v_1,v_2)_V = {}^tv_1 \bp &&1 \\ &1& \\ 1&& \ep \overline{v_2}.
\]
Let $\U(W)$ and $\U(V)$ be the isometry groups of $(W,(\cdot,\cdot)_W)$ and $(V,(\cdot,\cdot)_V)$ respectively. 
In the notation of \S\,\ref{SS:basic setting}, we have $\U(W) = \U_2(\R)$ and $\U(V) = \U_3(\R)$.
The tensor product $W\otimes_\C V$ admits a symplectic form over $\R$ defined by
\[
{\rm Tr}_{\C/\R}((\cdot,\cdot)_W\otimes (\cdot,\cdot)_V).
\]
Let $\psi_\R$ be the non-trivial additive character of $\R$ defined by $\psi_\R(x) = e^{2\pi\sqrt{-1}\,x}$ and $\psi_\C = \psi_\R\circ{\rm Tr}_{\C/\R}$. The additive character and symplectic form determine a Weil representation $\omega_{\psi_\R}$ of the metaplectic $\C^1$-cover ${\rm Mp}(W\otimes_\C V)$ of ${\rm Sp}(W\otimes_\C V)$. 
Let $\chi_V$  and $\chi_W$ be characters of $\C^\times$ such that 
\begin{align}\label{E:chi_V}
\chi_V\vert_{\R^\times} = {\rm sgn},\quad \chi_W\vert_{\R^\times} = 1.
\end{align}
Let $m_0$ and $n_0$ be the integers such that 
\[
\chi_V(z) = (z/\overline{z})^{m_0/2},\quad \chi_W(z) = (z/\overline{z})^{n_0/2}.
\]
Since $n_0$ is even, we can lift $\chi_W$ to a character $\tilde{\chi}_W$ of $\C^1$ defined by
\[
\tilde{\chi}_W(z/\overline{z}) = \chi_W(z),\quad z \in \C^\times.
\]
By the result of Kudla \cite{Kudla1994}, the data $(\psi_\R,\chi_V,\chi_W)$ specifying a splitting map 
\[
\U(W) \times \U(V) \rightarrow {\rm Mp}(W\otimes_\C V)(\R).
\]
%More precisely, in \cite[Theorem 3.1]{Kudla1994}, we take $\xi = \chi_V$ (resp.\,$\xi = \chi_V|\mbox{ }|_F^{-1/2}$) when $m$ is even (resp.\,odd).
Let $\omega_{\psi_\R,\chi_V,\chi_W}$ be the pullback of $\omega_{\psi_\R}$ to $\U(W) \times \U(V)$ on the Schr\"odinger model $\mathcal{S}({\rm M}_{3,1}(\C))$ which is the space of Bruhat--Schwartz functions on ${\rm M}_{3,1}(\C)$. 
We will abbreviate $\omega_{\psi_\R,\chi_V,\chi_W}$ as $\omega$.
The explicit formulas for $\omega$ are recalled as follows:
\begin{itemize}
\item For $h \in \U(V)$, we have
\begin{align}\label{E:Weil 1}
\omega(1,h)\varphi(x) = {\tilde{\chi}_W(\det h)}\varphi(h^{-1}x).
\end{align}
\item For $a \in \C^\times$, we have
\begin{align}\label{E:Weil 2}
\omega\left({\rm diag}(a,\overline{a}^{-1}),1 \right)\varphi(x) = \chi_V(a)|a|_\C^{3/2}\varphi(ax).
\end{align}
\item For $b \in \R$, we have
\begin{align}\label{E:Weil 3}
\omega\left(\bp 1 & b \\  & 1\ep,1 \right)\varphi(x) = \psi_\R\left(b(x_1\overline{x_3}+x_2\overline{x_2}+x_3\overline{x_1})\right)\varphi(x).
\end{align}
\item We have
\begin{align}\label{E:Weil 4}
\omega\left(\bp & 1 \\ -1 & \ep,1 \right)\varphi(x) = \sqrt{-1}\cdot\int_{{\rm M}_{3,1}(\C)}\varphi(y)\psi_\C(x_1\overline{y_3}+x_2\overline{y_2}+x_3\overline{y_1})\,dy,
\end{align}
%\begin{align}\label{E:Weil 4}
%\omega\left(\bp & -1 \\ 1 & \ep,1 \right)\varphi(x) = -\sqrt{-1}\cdot\int_{{\rm M}_{3,1}(\C)}\varphi(y)\psi_\C\left(-(x_1\overline{y_3}+x_2\overline{y_2}+x_3\overline{y_1})\right)\,dy,
%\end{align}
where $dy = dy_1dy_2dy_3$ with $dy_i$ twice the Lebesgue measure on $\C$ (self-dual with respect to $\psi_\C$).
%and 
%\[
%\gamma_{V,\psi} = \omega_{E/F}(\det Q)\cdot \gamma_F(-\delta^2,\psi)^m\cdot \gamma_F(-1,\psi)^{-m}
%\]
%with $\gamma_F(\cdot,\psi)$ the Weil index (cf.\,\cite[Appendix]{Rao1993}, \cite[\S\,A.1]{Ichino2005}) and a trace zero element $\delta \in E^\times$.
\end{itemize}

\subsubsection{Change of polarizations}

Let 
\begin{align*}
\mathcal{S}({\rm M}_{3,1}(\C)) \longrightarrow \mathcal{S}(\C) \otimes \mathcal{S}(\C^2),\quad \varphi  \longmapsto \hat{\varphi}
\end{align*}
be the partial Fourier transform defined by
\[
\hat{\varphi}(x;y) = \int_{\C}\varphi\bp z \\ x \\ y_1 \ep \psi_\C(z\overline{y_2})\,dz,\quad x \in \C,\,y=(y_1,y_2) \in \C^2.
\]
%for $x \in V_1(F)$ and $y=(y_1,y_2) \in E^2$. 
Here $dz$ is twice the Lebesgue measure on $\C$.
Let $\hat{\omega}$ be the representation of $\U(W) \times \U(V)$ on $\mathcal{S}(\C) \otimes \mathcal{S}(\C^2)$ defined by
\[
\hat{\omega}(g,h)\hat{\varphi} = (\omega(g,h)\varphi)\,\hat{}.
\]
We listed below some formulas for $\hat{\omega}$. The proof is by straightforward computation and we leave the details to the readers.
\begin{lemma}\label{E:partial Fourier}
Let $\varphi \in \mathcal{S}({\rm M}_{3,1}(\C))$. We have the following formulas:
\begin{itemize}
\item For $g \in \U(W)$ and $\hat{\varphi} = \varphi_1 \otimes \varphi_2$, we have
\begin{align}\label{E:PF 1}
\hat{\omega}(g,1)\hat{\varphi}(x;y) = \omega(g,1)\varphi_1(x)\cdot\varphi_2(yg).
\end{align}
Here the Weil representation of $\U(W) \times \C^1$ on $\mathcal{S}(\C)$ is defined as in \S\,\ref{SS:Weil rep} with respect to $(\psi_\R,\chi_V,\chi_W)$ and the standard hermitian form $z_1\overline{z_2}$ on $\C$.
\item For $u \in \C$, we have
\begin{align}\label{E:PF 2}
\hat{\omega}\left(1,\bp 1 & u & -\tfrac{1}{2}u\overline{u} \\  & 1 & -\overline{u} \\  &  & 1 \ep\right)\hat{\varphi}(x;y) = \psi_\C(ux\overline{y_2}+\tfrac{1}{2}u\overline{u} y_1\overline{y_2})\hat{\varphi}(x+\overline{u}y_1;y).
\end{align}
\item For $v = -\overline{v} \in \C$, we have
\begin{align}\label{E:PF 3}
\hat{\omega}\left(1,\bp 1 &  & v \\  & 1 &  \\  &  & 1 \ep\right)\hat{\varphi}(x;y) = \psi_\C(vy_1\overline{y_2})\hat{\varphi}(x;y).
\end{align}
\item For $a \in \C^\times$ and $z \in \C^1$, we have
\begin{align}\label{E:PF 4}
\hat{\omega}\left(1,{\rm diag}(a,z,\overline{a}^{-1})\right)\hat{\varphi}(x;y) = {\tilde{\chi}_W(a^{-1}\overline{a}z^{-1})}|a|_\C\hat{\varphi}(z^{-1}x;\overline{a}y).
\end{align}
\end{itemize}
\end{lemma}

\subsubsection{Explicit theta correspondence}\label{SS:explicit theta}

For an irreducible representation $\pi$ of $\U(W)$, the maximal $\pi$-isotypic quotient of $\omega = \omega_{\psi_\R,\chi_V,\chi_W}$ is of the form
\[
\pi \boxtimes \Theta_{\psi_\R,\chi_V,\chi_W}(\pi)
\]
for some smooth representation $\Theta_{\psi_\R,\chi_V,\chi_W}(\pi)$ of $\U(V)$.
By the result of Howe \cite{Howe1989} (Howe duality principle), the maximal semisimple quotient $\theta_{\psi_\R,\chi_V,\chi_W}(\pi)$ of $\Theta_{\psi_\R,\chi_V,\chi_W}(\pi)$ is either zero or irreducible. We write $\theta(\pi) = \theta_{\psi_\R,\chi_V,\chi_W}(\pi)$ and call it the theta lift of $\pi$ (with respect to the data $(\psi_\R,\chi_V,\chi_W)$) whenever it is non-zero.

Let $N_2(\R)$ and $N_3(\R)$ be the maximal unipotent subgroups of $\U(W)$ and $\U(V)$ defined in \S\,\ref{SS:basic setting}. We then have the non-degenerated characters $\psi_{N_2}$ and $\psi_{N_3}$ of $N_2(\R)$ and $N_3(\R)$ respectively. Recall that
\begin{align*}
\psi_{N_2}\left( \bp 1 & x \\ & 1 \ep  \right) = \psi_\R(x),\quad \psi_{N_3} \left( \bp 1 & x & y-\tfrac{1}{2}x\overline{x} \\  & 1 & -\overline{x} \\  &  & 1 \ep \right) = \psi_\C(x).
\end{align*}
Let $C^\infty(N_3(\R)\backslash \U(V),\psi_{N_3})^{\rm mg}$ be the space of smooth functions $W: \U(V) \rightarrow \C$ which are of moderate growth and satifying $W(nh) = \psi_{N_3}(n)W(h)$ for $n \in N_3(\R)$ and $h \in \U(V)$.

Let $\pi$ be an irreducible $\psi_{N_2}$-generic representation of $\U(W)$. Then $\pi^\vee$ is $\psi_{N_2}^{-1}$-generic. Recall $\mathcal{W}(\pi^\vee,\psi_{N_2}^{-1})$ is the space of Whittaker functions of $\pi^\vee$ with respect to $\psi_{N_2}^{-1}$. 
Consider an equivariant homomorphism
\begin{align}\label{E:explicit theta}
\begin{split}
\mathcal{L}_\pi:\mathcal{W}(\pi^\vee,\psi_{N_2}^{-1}) \times \mathcal{S}({\rm M}_{3,1}(\C)) \longrightarrow C^\infty(N_3(\R)\backslash \U(V),\psi_{N_3})^{\rm mg},\\
\mathcal{L}_\pi(h;W,\varphi):= \int_{N_2(\R)\backslash \U(W)} \hat{\omega}(g,h)\hat{\varphi}(1;0,1)W(g)\,\overline{dg},\quad h\in \U(V).
\end{split}
\end{align}
Here $\overline{dg} = dx \backslash dg$ is the quotient measure on $N_2(\R)\backslash \U(W)$, with 
\begin{align}\label{E:U(W) measure}
dg = y^{-3}dx\,dy\,dk
\end{align}
for 
\[
g = \bp 1 & x \\ & 1 \ep {\rm diag}(y,y^{-1})k,\quad x \in \R,\, y \in \R_+,\, k \in \U(2)\cap\U(W),
\]
$dx$, $dy$ are the Lebesgue measure and ${\rm vol}(\U(2)\cap\U(W),dk)=1$.
%\begin{align}\label{E:U(W) measure}
%dg = |a|_\C^{-1}dx\,d^\times a\,dk
%\end{align}
%for 
%\[
%g = \bp 1 & x \\ & 1 \ep {\rm diag}(a,\overline{a}^{-1})k,\quad x \in \R,\, a \in \C^\times,\, k \in \U(2)\cap\U(W),
%\]
%$dx$ is the Lebesgue measure, ${\rm vol}(\U(2)\cap\U(W),dk)=1$, and $d^\times a = (2\pi)^{-1}r^{-1}drd\theta$
%for $a = r e^{\sqrt{-1}\,\theta}$ with $r \in \R_+$ and $\theta \in [0,2\pi]$.
The equivariance under $\psi_{N_3}$ follows immediately from (\ref{E:PF 2}) and (\ref{E:PF 3}).
For the absolute convergence and growth condition, we include a proof in the following lemma. 
\begin{lemma}\label{L:explicit theta}
Let $\pi$ be an irreducible $\psi_{N_2}$-generic representation of $\U(W)$.
\begin{itemize}
\item[(1)] The integrations defining $\mathcal{L}_\pi$ converge absolutely and are of moderate growth.
\item[(2)] If $\pi$ is a discrete series representation and $\mathcal{L}_\pi$ is non-zero, then the image of $\mathcal{L}_\pi$ is equal to $\mathcal{W}(\theta(\pi),\psi_{N_3})$.
\end{itemize}
\end{lemma}

\begin{proof}
Let $W \in \mathcal{W}(\pi^\vee,\psi_{N_2}^{-1})$ and $\varphi \in \mathcal{S}({\rm M}_{3,1}(\C))$. Let $C,M$ be constants such that
\[
|W({\rm diag}(y,y^{-1})k)| \leq C y^M,\quad y \in \R_+,\, k \in \U(2)\cap \U(W).
\]
Also there exists a positive Bruhat--Schwartz function $\varphi_1 \otimes \varphi_2 \in \mathcal{S}(\C) \otimes \mathcal{S}(\C^2)$ such that
\[
|\hat{\omega}(k_1,k_2)\hat{\varphi}| \leq \varphi_1\otimes\varphi_2,\quad (k_1,k_2) \in (\U(2)\cap \U(W)) \times (\U(3) \cap \U(V)).
\]
By (\ref{E:PF 1}) and (\ref{E:PF 4}), we have
\[
|\hat{\omega}\left( {\rm diag}(y,y^{-1}),{\rm diag}(a^{-1},1,a)\right)(\varphi_1\otimes\varphi_2)(1;0,1)| = a^2y\varphi_1(y)\varphi_2(0,ay^{-1}),\quad a,y \in \R_+.
\]
For $a,y \in \R_+$ and $k \in \U(3)\cap\U(V)$, we thus have
\begin{align*}
&|\mathcal{L}_\pi({\rm diag}(a,1,a^{-1})k;W,\varphi)|\\
&\leq C\int_0^\infty \left|\hat{\omega}\left( {\rm diag}(y,y^{-1}),{\rm diag}(a,1,a^{-1})\right)(\varphi_1\otimes\varphi_2)(1;0,1)\right|y^{M-2}\,d y\\
& = Ca^2\int_0^\infty \varphi_1(y)\varphi_2(0,ay^{-1})y^{M-2}\,d y.
\end{align*}
Let $C_1$ and $C_2$ be upper bounds of $\varphi_1$ and $\varphi_2$ respectively. 
Fix positive integers $n_1$ and $n_2$ such that 
\[
n_1>M-1,\quad n_2>-M+2.
\]
Let $D_1$ and $D_2$ be constants such that 
\[
\varphi_1(x) \leq D_1|x|^{-n_1},\quad \varphi_2(0,x) \leq D_2 |x|^{-n_2},\quad x \in \R.
\]
Then
\begin{align*}
\int_a^\infty \varphi_1(y)\varphi_2(0,ay^{-1})y^{M-2}\,d y &\leq C_2D_1\int_a^\infty y^{-n_1+M-2}\,dy = C_2D_1 \frac{a^{-n_1+M-1}}{n_1-M+1}
\end{align*}
and 
\begin{align*}
\int_0^a \varphi_1(y)\varphi_2(0,ay^{-1})y^{M-2}\,d y \leq C_1D_2 a^{-n_2} \int_0^a y^{n_2+M-2}\,d y = C_1D_2\frac{a^{M-1}}{n_2+M-1}.
\end{align*}
In conclusion, we have
\[
|\mathcal{L}_\pi({\rm diag}(a,1,a^{-1})k;W,\varphi)| \leq C\min\left\{C_2D_1 \frac{a^{-n_1+M+1}}{n_1-M+2},C_1D_2\frac{a^{M+1}}{n_2+M-1}\right\}.
\]

For the second assertion, by the same arguments in the proof of \cite[Proposition C.1-(i)]{GI2016}, the condition that $\pi$ is a discrete series representation implies that the big theta lift $\Theta(\pi)$ is semisimple. Therefore, $\Theta(\pi) = \theta(\pi)$.
If $\mathcal{L}_\pi$ is non-zero, then it induces a surjective non-zero intertwining map from $\Theta(\pi)$ to ${\rm image}(\mathcal{L}_\pi)$, which must be an isomorphism by the irreducibility of $\Theta(\pi)$.
This completes the proof.
\end{proof}

Define $k_W \in {\rm SU}(2)$ and $k_V \in {\rm SU}(3)$ by
\[
k_W = \bp \tfrac{1}{\sqrt{2}} & \tfrac{-\sqrt{-1}}{\sqrt{2}} \\ \tfrac{-\sqrt{-1}}{\sqrt{2}} & \tfrac{1}{\sqrt{2}} \ep,\quad k_V = \bp \tfrac{1}{\sqrt{2}} & 0 &\tfrac{1}{\sqrt{2}}\\0&1&0\\\tfrac{-1}{\sqrt{2}} & 0 & \tfrac{1}{\sqrt{2}} \ep.
\]
Then we have isomorphisms
\begin{align}\label{E:iso. 2}
\begin{split}
\U(W) &\longrightarrow \U(1,1),\quad g \longmapsto k_W g k_W^{-1}\\
\U(V) &\longrightarrow \U(2,1),\quad  g\longmapsto k_V g k_V^{-1}.
\end{split}
\end{align}
Under these isomorphisms, we parameterize discrete series representations of $\U(W)$ and $\U(V)$ by that of $\U(1,1)$ and $\U(2,1)$ as in \S\,\ref{SS:d.s. rep.}.
Let $\pi = \pi_\mu$ be a discrete series representation of $\U(W)$ with Harish-Chandra parameter $\mu = (\mu_1,\mu_2) \in (\Z+\tfrac{1}{2})^2$. Note that $\pi_\mu$ is $\psi_{N_2}$-generic if and only if $\mu_1 > \mu_2$ (cf.\,\cite[Proposition A.3]{Atobe2020}).
By the result of Paul \cite{Paul2000} on local theta correspondence for real unitary groups in the almost equal rank case (see also \cite{Ichino2022}), the theta lift $\theta(\pi)$ is a discrete series representation if and only if $\mu_1+\tfrac{m_0}{2} \neq 0$ and $\mu_2+\tfrac{m_0}{2} \neq 0$. 
More explicitly, if $\mu_1>\mu_2$, $\mu_1+\tfrac{m_0}{2} \neq 0$, and $\mu_2+\tfrac{m_0}{2} \neq 0$, then we have 
\[
\theta(\pi_\mu) = \pi_{\lambda(\mu)},
\]
where
\[
\lambda(\mu) = \begin{cases}
(\tfrac{-n_0}{2},\mu_2+\tfrac{m_0-n_0}{2},\mu_1+\tfrac{m_0-n_0}{2}) & \mbox{ if $\mu_1+\tfrac{m_0}{2}<0$},\\
(\mu_1+\tfrac{m_0-n_0}{2},\mu_2+\tfrac{m_0-n_0}{2},\tfrac{-n_0}{2}) & \mbox{ if $\mu_2+\tfrac{m_0}{2} < 0 <\mu_1+\tfrac{m_0}{2}$},\\
(\mu_1+\tfrac{m_0-n_0}{2},\tfrac{-n_0}{2},\mu_2+\tfrac{m_0-n_0}{2}) & \mbox{ if $0<\mu_2+\tfrac{m_0}{2}$}.
\end{cases}
\]
In particular, to obtain all generic discrete series representations of $\U(V)$ from the theta correspondence, it is sufficient to consider the first or third types.
From now on, we assume 
\[
\mu_1>\mu_2,\quad \mu_1+\tfrac{m_0}{2} < 0.
\]
Let 
\begin{align}\label{E:theta parameter}
(\alpha_1,\alpha_3,\alpha_2):=\lambda(\mu) = (\tfrac{-n_0}{2},\mu_2+\tfrac{m_0-n_0}{2},\mu_1+\tfrac{m_0-n_0}{2}),
\end{align}
which is compatible with the notation for generic discrete series representations in (\ref{E:generic parameter}).
Define a Bruhat--Schwartz function $\varphi_\mu \in \mathcal{S}({\rm M}_{3,1}(\C))$ by 
\begin{align}\label{E:archimedean test function}
\varphi_\mu(x) = (x_1-x_3)^{\alpha_1-\alpha_2}x_2^{\alpha_1-\alpha_3}e^{-2\pi {}^t\!x^cx}.
\end{align}
Note that both $\lambda(\mu)$ and $\varphi_\mu$ depend on the choice of $m_0$ and $n_0$ (that is, $\chi_V$ and $\chi_W$), as the theta correspondence depending on these data.
It is not difficult to verify that $\varphi_\mu$ is a highest weight vector under the action of $(\U(2)\cap\U(W)) \times (\U(3)\cap\U(V))$. More precisely, we have the following lemma which can be verified directly by the explicit formulas of the Weil representations and we leave the details to the readers.

\begin{lemma}\label{L:archimedean test function}
\noindent
\begin{itemize}
\item[(1)] 
For $\theta_1,\theta_2,\theta_3 \in \R$, we have 
\begin{align*}
\omega \left(k_W^{-1}{\rm diag}(e^{\sqrt{-1}\theta_1},e^{\sqrt{-1}\theta_2})k_W,1\right)\varphi_\mu &=  e^{\sqrt{-1}\,((-\mu_2+1/2)\theta_1+ (-\mu_1-1/2) \theta_2)}\varphi_\mu,\\
\omega\left(1,\,k_V^{-1}{\rm diag}(e^{\sqrt{-1}\theta_1},e^{\sqrt{-1}\theta_2},e^{\sqrt{-1}\theta_3})k_V\right)\varphi_\mu &= e^{\sqrt{-1}\,(\alpha_1\theta_1+\alpha_3\theta_2+\alpha_2\theta_3)}\varphi_\mu.
\end{align*}
\item[(2)] Let
\[
N_+ = \bp 0 & 1 & 0 \\ -1 & 0 & 0 \\ 0&0&0 \ep\otimes\frac{1}{2} + \bp 0 & -\sqrt{-1} & 0 \\ -\sqrt{-1} & 0 & 0 \\ 0&0&0\ep \otimes \frac{\sqrt{-1}}{2} \in {\rm Lie}(\U(2) \times \U(1))_\C.
\]
We have
\[
\omega\left( 1,k_V^{-1}N_+ k_V \right)\varphi_\mu=0.
\]
%\item[(3)] For $\theta \in \R$, we have
%\begin{align*}
%\omega_{\psi_\infty,\chi_{V,\infty}}\left( \bp \cos\theta & \sin\theta\\-\sin\theta & \cos\theta\ep ,1\right)\varphi_\infty^\iota &= e^{\sqrt{-1}(-\lambda_2+\lambda_3+1)\theta}\varphi_\infty^\iota,\\
%\omega_{\psi_\infty,\chi_{V,\infty}}\left(e^{\sqrt{-1}\theta},1\right)\varphi_\infty^\iota &= e^{\sqrt{-1}(-\lambda_2-\lambda_3+2m_0)\theta}\varphi_\infty^\iota.
%\end{align*}
\end{itemize}
\end{lemma}

In the following lemma, we show that $\mathcal{L}_{\pi_\mu}$ is non-zero. In particular, by Lemma \ref{L:explicit theta}-(2), it defines an intertwining map from $\pi_{\mu}^\vee \otimes \omega$ to $\pi_{\lambda(\mu)}$.
More precisely, let $W_{\pi_\mu^\vee} \in \mathcal{W}(\pi_\mu^\vee,\psi_{N_2}^{-1})$ be a lowest weight Whittaker function in the minimal $(\U(2) \cap \U(W))$-type of $\pi_\mu^\vee$ normalized so that $W_{\pi_\mu^\vee}(1) = e^{-2\pi}$. The explicit formula is given by
\begin{align}\label{E:archimedean test function 2}
W_{\pi_\mu^\vee}\left( {\rm diag}(a,a^{-1})k_W^{-1}{\rm diag}(e^{\sqrt{-1}\theta_1},e^{\sqrt{-1}\theta_2})k_W  \right) = e^{\sqrt{-1}\,((\mu_2-1/2)\theta_1+ (\mu_1+1/2) \theta_2)}a^{\mu_1-\mu_2+1}e^{-2\pi a^2},
\end{align}
for $a>0$ and $\theta_1,\theta_2 \in \R$.
By Lemma \ref{L:archimedean test function}, the Whittaker function $\mathcal{L}_{\pi_\mu}(W_{\pi_\mu^\vee},\varphi_\mu)$, provided it is non-zero, defines a highest weight Whittaker function in the minimal $(\U(3) \cap \U(V))$-type of $\pi_{\lambda(\mu)}$.

\begin{lemma}\label{L:local Whittaker 2}
For $a>0$, we have
\begin{align*}
\mathcal{L}_{\pi_\mu}\left({\rm diag}(a,1,a^{-1});W_{\pi_\mu^\vee},\varphi_\mu\right) &= (\sqrt{-1})^{\alpha_1-\alpha_2}(\sqrt{2})^{-\alpha_2+\alpha_3}\cdot a^{\alpha_1-\alpha_3+2}K_{\alpha_2-\alpha_3}(4\sqrt{2}\pi a).
\end{align*}
\end{lemma}

\begin{proof}
It is easy to verify that
\[
\hat{\varphi}_\mu(x;y) = x^{\alpha_1-\alpha_3}(\sqrt{-1}y_2-y_1)^{\alpha_1-\alpha_2}e^{-2\pi(xx^c+y{}^ty^c)}.
\]
For $y >0$, by (\ref{E:PF 1}) and (\ref{E:PF 4}) we have
\begin{align*}
\hat{\omega}\left( {\rm diag}(y,y^{-1}),{\rm diag}(a,1,a^{-1})\right)\hat{\varphi}_\mu(1;0,1)& = a^2y\hat{\varphi}_\mu(y;0,ay^{-1})\\
& = (\sqrt{-1})^{\alpha_1-\alpha_2}\cdot a^{\alpha_1-\alpha_2+2}y^{\alpha_2-\alpha_3+1}e^{-2\pi(y^2+a^2y^{-2})}.
\end{align*}
%Note that
%\[
%-\lambda_2+\lambda_3+1 = \kappa,\quad -\lambda_2-\lambda_3+2m_0=-{\sf w}.
%\]
We thus conclude from Lemma \ref{L:archimedean test function}-(1) and (\ref{E:archimedean test function 2}) that
\begin{align*}
&\mathcal{L}_{\pi_\mu}\left({\rm diag}(a,1,a^{-1});W_{\pi_\mu^\vee},\varphi_\mu\right)\\
& = \int_0^\infty \hat{\omega}\left( {\rm diag}(y,y^{-1}),{\rm diag}(a,1,a^{-1})\right)\hat{\varphi}_\mu(1;0,1)W_{\pi_\mu^\vee}\left( {\rm diag}(y,y^{-1}) \right)y^{-3}\,dy\\
&= (\sqrt{-1})^{\alpha_1-\alpha_2}\cdot a^{\alpha_1-\alpha_2+2}\int_0^\infty y^{2\alpha_2-2\alpha_3-1}e^{-2\pi(2y^2+a^2y^{-2})}\,d y\\
%& = 2^{-1}a^{\alpha_1-\lambda_3+2}\int_{0}^\infty y^{-\lambda_2+\lambda_3} e^{-4\sqrt{2}\pi a(y+y^{-1})/2}\,d^\times y\\
&= (\sqrt{-1})^{\alpha_1-\alpha_2}(\sqrt{2})^{-\alpha_2+\alpha_3}\cdot a^{\alpha_1-\alpha_3+2}K_{\lambda_2-\lambda_3}(4\sqrt{2}\pi a).
\end{align*}
%Note that $\overline{\omega}_{\psi_\infty,\chi_{V,\infty}} = \omega_{\psi_\infty^{-1},\chi_{V,\infty}^c}$. Hence the second formula follows from the first one and Lemma \ref{L:archimedean test function}-(1).
This completes the proof.
\end{proof}

\subsubsection{Doubling local zeta integrals}

In this section, we compute certain local zeta integrals which appear in \cite[Proposition 4.1]{Morimoto2024}. These integrals are (archimedean) local part of the Rallis inner product formula which relates the Petersson norm of global theta lifts with special values of standard $L$-functions on unitary groups.

Let $\pi_\mu$ be the discrete series representation of $\U(W)$ with Harish-Chandra parameter $\mu \in (\Z+\tfrac{1}{2})^2$. 
As in \S\,\ref{SS:explicit theta}, we assume $\mu_1>\mu_2$ and $\mu_1+\tfrac{m_0}{2}<0$.
We fix an equivariant hermitian pairing $(\cdot,\cdot)$ on $\pi_\mu^\vee$.
Let $\mathcal{B} : \mathcal{S}({\rm M}_{3,1}(\C)) \times \mathcal{S}({\rm M}_{3,1}(\C)) \rightarrow \C$ be the hermitian pairing defined by
\[
\mathcal{B} (\varphi_1,\varphi_2) = \int_{\C^3}\varphi_1(z)\overline{\varphi_2(z)}\,dz,
\]
where $dz = dz_1\,dz_2\,dz_3$ with $dz_i$ twice the Lebesgue measure on $\C$ (self-dual with respect to $\psi_\C$).
Note that $\mathcal{B}$ is equivariant under the Weil representation $\omega \otimes \overline{\omega}$.
For $f_1,f_2 \in \pi_\mu^\vee$ and $\varphi_1,\varphi_2 \in \mathcal{S}({\rm M}_{3,1}(\C))$, we define the doubling local zeta integral
\begin{align*}
Z^*(f_1,f_2,\varphi_1,\varphi_2) = \frac{\Gamma_\R(2)\Gamma_\R(4)}{\Gamma_\C(1+\alpha_1-\alpha_2)\Gamma_\C(1+\alpha_1-\alpha_3)}\int_{\U(W)}\mathcal{B}( \omega(g,1)\varphi_1,\varphi_2)(\pi_\mu^\vee(g)f_1,f_2)\,dg.
\end{align*}
Here $dg$ is the Haar measure defined in (\ref{E:U(W) measure}), and $(\alpha_1,\alpha_2,\alpha_3) \in \Z^3$ is defined in (\ref{E:theta parameter}).
Note that the product of gamma factors in the denominator is the value at $s=1$ of the $\chi_V$-twisted standard $L$-function of $\pi_\mu$:
\[
L(s,\pi_\mu,{\rm std}\otimes\chi_V) = \Gamma_\C(s+\alpha_1-\alpha_2)\Gamma_\C(s+\alpha_1-\alpha_3).
\]

\begin{lemma}\label{L:local zeta 2}
Let $f_{\pi_\mu^\vee} \in \pi_\mu^\vee$ be a lowest weight vector in the minimal $(\U(2)\cap \U(W))$-type of $\pi_\mu^\vee$. We have
\begin{align*}
\frac{Z^*(f_{\pi_\mu^\vee},f_{\pi_\mu^\vee},\varphi_\mu,\varphi_\mu)}{(f_{\pi_\mu^\vee},f_{\pi_\mu^\vee})}
 = 2^{-\alpha_1+\alpha_3-1}(\alpha_1-\alpha_3+1)^{-1}.
\end{align*}
\end{lemma}

\begin{proof}
By Lemma \ref{L:archimedean test function}-(1) and the Cartan decomposition, we have
\begin{align*}
&\int_{\U(W)}\mathcal{B}(\omega(g,1)\varphi_\mu,\varphi_\mu)(\pi_\mu^\vee(g)f_{\pi_\mu^\vee},f_{\pi_\mu^\vee})\,dg \\ 
&= {2\pi}\int_0^\infty \mathcal{B}(\omega({\rm diag}(e^t,e^{-t}),1)\varphi_\mu,\varphi_\mu)(\pi_\mu^\vee({\rm diag}(e^t,e^{-t}))f_{\pi_\mu^\vee},f_{\pi_\mu^\vee})\sinh(2t)\,dt.
\end{align*}
For $a>0$, we have
\begin{align*}
&\mathcal{B}(\omega({\rm diag}(a,a^{-1}),1)\varphi_\mu,\varphi_\mu) \\ 
& = a^{2\alpha_1-\alpha_2-\alpha_3+3}\int_{\C^3}|x_1-x_3|_\C^{\alpha_1-\alpha_2}|x_2|_\C^{\alpha_1-\alpha_3}e^{-2\pi {}^tx^cx (a^2+1)}\,dx\\
& = 2^{\alpha_1-\alpha_2+1}(a+a^{-1})^{-2\alpha_1+\alpha_2+\alpha_3-3} \int_{\C^3}|y_2|_\C^{\alpha_1-\alpha_2}|x_2|_\C^{\alpha_1-\alpha_3}e^{-2\pi(|y_1|_\C+|y_2|_\C+|x_2|_\C)}\,dy_1\,dy_2\,dx_2\\
& = 2^{\alpha_1-\alpha_2-1}(a+a^{-1})^{-2\alpha_1+\alpha_2+\alpha_3-3} (2\pi)^2\Gamma_\C(1+\alpha_1-\alpha_2)\Gamma_\C(1+\alpha_1-\alpha_3).
\end{align*}
Here the second equality follows from the change of variables $x_1+x_3 = \sqrt{2}y_1$ and $x_1-x_3 = \sqrt{2}y_2$.
By the formula of matrix coefficients for $\SL_2(\R)$, we have
\[
\frac{(\pi_\mu^\vee({\rm diag}(a,a^{-1}))f_{\pi_\mu^\vee},f_{\pi_\mu^\vee})}{(f_{\pi_\mu^\vee},f_{\pi_\mu^\vee})} = \left( \frac{a+a^{-1}}{2}\right)^{-\alpha_2+\alpha_3-1},\quad a>0.
\]
Therefore, we have
\begin{align*}
&\mathcal{B}(\omega({\rm diag}(a,a^{-1}),1)\varphi_\mu,\varphi_\mu) \\ 
& = 2^{-\alpha_1+\alpha_3-1}\frac{\Gamma_\C(1+\alpha_1-\alpha_2)\Gamma_\C(1+\alpha_1-\alpha_3)}{\Gamma_\R(2)\Gamma_\R(4)} \int_0^\infty\cosh(t)^{-2\alpha_1+2\alpha_3-4}\sinh(2t)\,dt\\
& = 2^{-\alpha_1+\alpha_3-1}\frac{\Gamma_\C(1+\alpha_1-\alpha_2)\Gamma_\C(1+\alpha_1-\alpha_3)}{\Gamma_\R(2)\Gamma_\R(4)}(\alpha_1-\alpha_3+1)^{-1}.
\end{align*}
This complete the proof.

\end{proof}

\subsubsection{Computation of $C_{\itPi_\infty,\psi_{N,\infty}}$}

Recall the global unitary group $G$ over $\Q$ defined in \S\,\ref{E:basic}.
We have identified $G(\R)$ with $\U(2,1)$ and $\U(V)$ with $\U(2,1)$ by the isomorphisms in (\ref{E:iso. 1}) and (\ref{E:iso. 2}) respectively. These induce an isomorphism from $G(\R)$ to $\U(V)$ which is given by
\[
G(\R) \longrightarrow \U(V),\quad g \longmapsto {\rm diag}(-\delta^{-1},1,1)g{\rm diag}(-\delta,1,1).
\]
Under the isomorphisms $G(\R) \cong \U(V) \cong \U(2,1)$, we have
\[
K_\infty \cong \U(3)\cap \U(V) \cong \{{\rm diag}(k_1,k_2)\,\vert\,k_1 \in \U(2),\,k_2\in\U(1)\}.
\]

Recall the archimedean component $\itPi_\infty = \pi_\lambda$ for some $\lambda = (\alpha_1,\alpha_3,\alpha_2)$ with $\alpha_1 > \alpha_2 > \alpha_3$. We choose $\chi_V$, $\chi_W$, and $\mu = (\mu_1,\mu_2) \in (\Z+\tfrac{1}{2})^2$ such that
\[
\mu_1>\mu_2,\quad \mu_1+\tfrac{m_0}{2}<0, \quad \lambda(\mu) = \lambda.
\]
Then we have 
\[
\theta(\pi_\mu) = \pi_\lambda
\]
by the result of Paul \cite{Paul2000}. By Lemmas \ref{L:explicit theta} and \ref{L:local Whittaker 2}, we have an explicit realization $\mathcal{L}_{\pi_\mu}$ (cf.\,(\ref{E:explicit theta})) of the theta lifting.
Fix equivariant hermitian pairings $\<\cdot,\cdot\>$ and $(\cdot,\cdot)$ on $\pi_\lambda$ and $\pi_\mu^\vee$ respectively. 
Let $I_\lambda$ (resp.\,$I_{\mu^\vee}$) be the Whittaker functional on $\pi_\lambda \otimes \overline{\pi_\lambda}$ (resp.\,$\pi_\mu^\vee \otimes \overline{\pi_\mu^\vee}$) defined as in (\ref{E:local Whittaker 1}) with respect to $\<\cdot,\cdot\>$ and $\psi_{N_3}$ (resp.\,$(\cdot,\cdot)$ and $\psi_{N_2}$), where the Haar measures on $N_3(\R)$ and $N_2(\R)$ are normalized as follows. 
\begin{itemize}
\item For $N_3(\R)$, we have $dn = dxdy$ for $n=n(x,y)$ in (\ref{E:additive character}), where $dx$ is twice the Lebesgue measure on $\C$ and $dy$ is the pushforward of the Lebesgue measure on $\R$ under the isomorphism $y \mapsto \sqrt{-1}y$ from $\R$ to $\sqrt{-1}\cdot\R$.
\item For $N_2(\R)$, $dn=dx$ is the Lebesgue measure on $\R$ for $n = \bp 1 & x \\ & 1\ep$.
\end{itemize}
We define the normalized functionals
\[
I_\lambda^* := \frac{L(1,\pi_\lambda,{\rm Ad})}{\Gamma_\R(2)^2\Gamma_\R(4)}\cdot I_\lambda,\quad I_{\mu^\vee}^* := \frac{L(1,\pi_\mu^\vee,{\rm Ad})}{\Gamma_\R(2)^2}\cdot I_{\mu^\vee}.
\]
The local constants we are interested in are: 
\begin{align}\label{E:archimedean constants}
C_{\pi_\lambda,\psi_{N_3}} := \frac{I_\lambda^*(f_{\pi_\lambda}\otimes f_{\pi_\lambda})}{\<f_{\pi_\lambda},f_{\pi_\lambda}\>},\quad C_{\pi_{\mu^\vee},\psi_{N_2}} := \frac{I_{\mu^\vee}^*(f_{\pi_{\mu^\vee}}\otimes f_{\pi_{\mu^\vee}})}{(f_{\pi_{\mu^\vee}},f_{\pi_{\mu^\vee}})}.
\end{align}
Here $f_{\pi_\lambda} \in \pi_\lambda$ (resp.\,$f_{\pi_{\mu^\vee}} \in \pi_{\mu^\vee}$) is a highest weight (resp.\,lowest weight) vector in the minimal $(\U(3)\cap \U(V))$-type (resp.\,$(\U(2)\cap \U(W))$-type) of $\pi_\lambda$ (resp.\,$\pi_{\mu^\vee}$).

\begin{lemma}\label{L:change of character}
Let $C_{\itPi_\infty,\psi_{N,\infty}}$ be the archimedean local constant in (\ref{E:local constants}).
We have
\[
C_{\itPi_\infty,\psi_{N,\infty}} = D_\K^2\cdot C_{\pi_\lambda,\psi_{N_3}}.
\]
\end{lemma}

\begin{proof}
Note that
\[
\psi_{N,\infty}(n(\delta x, \delta^2 y)) = \psi_\C(x) = \psi_{N_3}({\rm diag}(-\delta^{-1},1,1)n(x,y){\rm diag}(-\delta,1,1)).
\]
The assertion then follows from the change of variables 
\[
n(x,y) \longmapsto {\rm diag}(\delta,1,\delta^{-1})n(x,y){\rm diag}(\delta^{-1},1,\delta) = n(\delta x, \delta^2 y).
\]
\end{proof}
 
To compute $C_{\pi_\lambda,\psi_{N_3}}$, we use the result of Morimoto \cite[Proposition 4.1]{Morimoto2024} to express it in terms of $C_{\pi_{\mu^\vee},\psi_{N_2}}$, $\mathcal{L}_{\pi_\mu}$, and the doubling zeta integral $Z^*$.
More precisely, in $loc.$ $cit.$, the following is proved:
Fix a non-zero equivariant map
\[
\theta : \pi_\mu^\vee \otimes \omega \longrightarrow \pi_\lambda.
\]
Let $f \in \pi_{\mu^\vee}$ and $\varphi \in \mathcal{S}({\rm M}_{3,1}(\C))$. Let $W_f \in \mathcal{W}(\pi_\mu^\vee,\psi_{N_2}^{-1})$ be a Whittaker function associated with $f$. Assume $W_f(1) \neq 0$ and $\theta(f,\varphi) \neq 0$. Then we have
\begin{align}\label{E:Morimoto}
\frac{Z^*(f,f,\varphi,\varphi)}{(f,f)}\cdot \frac{I_\lambda^*(\theta(f,\varphi) \otimes \theta(f,\varphi))}{\<\theta(f,\varphi),\theta(f,\varphi)\>} = \left\vert\frac{\mathcal{L}_{\pi_\mu}(W_f,\varphi)}{W_f(1)}\right\vert^2\cdot \frac{I_{\mu^\vee}^*(f\otimes f)}{(f,f)}.
\end{align}

\begin{prop}\label{P:archimedean}
We have
\[
C_{\itPi_\infty,\psi_{N,\infty}} = D_\K^2\cdot 2^{\alpha_1-2\alpha_2+\alpha_3-1}(\alpha_1-\alpha_3+1)K_{\alpha_2-\alpha_3}(4\sqrt{2}\pi)^2.
\]
\end{prop}

\begin{proof}
Let $\varphi_\mu$ be the Bruhat--Schwartz function defined in (\ref{E:archimedean test function}).
By Lemma \ref{L:archimedean test function}, we have
\[
f_{\pi_\lambda}:=\theta(f_{\pi_\mu^\vee},\varphi_\mu)
\]
is a highest weight vector in the minimal $(\U(3)\cap \U(V))$-type of $\pi_\lambda$.
Let $W_{\pi_{\mu^\vee}}$ be the Whittaker function associated with $f_{\pi_{\mu^\vee}}$ normalized as in (\ref{E:archimedean test function 2}).
Apply (\ref{E:Morimoto}) to $f=f_{\pi_{\mu^\vee}}$ and $\varphi = \varphi_\mu$, by Lemmas \ref{L:local Whittaker 2}, \ref{L:local zeta 2}, \ref{L:change of character}, we then have
\[
2^{-\alpha_1+\alpha_3-1}(\alpha_1-\alpha_3+1)^{-1}\cdot D_\K^{-2}\cdot C_{\itPi_\infty,\psi_{N,\infty}} = \frac{2^{-\alpha_2+\alpha_3}K_{\alpha_2-\alpha_3}(4\sqrt{2}\pi)^2}{e^{-4\pi}}\cdot C_{\pi_{\mu^\vee},\psi_{N_2}}.
\]
It thus remains to show that
\[
C_{\pi_{\mu^\vee},\psi_{N_2}} = 2^{-\alpha_2+\alpha_3-2}e^{-4\pi}.
\]
We normalized the matrix coefficient $(\pi_{\mu^\vee}(g)f_{\pi_{\mu^\vee}},f_{\pi_{\mu^\vee}})$ so that it is given by %(cf.\,{\color{red}??})
\[
(\pi_{\mu^\vee}(g)f_{\pi_{\mu^\vee}},f_{\pi_{\mu^\vee}}) = \int_0^\infty W_{\pi_{\mu^\vee}}\left({\rm diag}\left(\sqrt{a},\sqrt{a}^{-1}\right)g\right) \overline{W_{\pi_{\mu^\vee}}\left({\rm diag}\left(\sqrt{a},\sqrt{a}^{-1}\right) \right)}\,\frac{da}{a},\quad g \in \U(W).
\]
In particular, 
\begin{align*}
\left(\pi_{\mu^\vee}\left( \bp 1 & x \\ & 1 \ep\right)f_{\pi_{\mu^\vee}},f_{\pi_{\mu^\vee}}\right) 
& = \int_0^\infty a^{\alpha_2-\alpha_3+1}e^{-4\pi a+2\pi\sqrt{-1}\,ax}\,\frac{da}{a}\\
& = (4\pi-2\pi\sqrt{-1}\,x)^{-\alpha_2+\alpha_3-1}\Gamma(\alpha_2-\alpha_3+1).
\end{align*}
Therefore, by the residue theorem we have
\begin{align*}
I_{\mu^\vee}(f_{\pi_{\mu^\vee}} \otimes f_{\pi_{\mu^\vee}}) &= 2 (-4\pi\sqrt{-1})^{-\alpha_2+\alpha_3-1}\Gamma(\alpha_2-\alpha_3+1) \int_\R \frac{e^{-4\pi\sqrt{-1}\,x}}{(x+\sqrt{-1})^{\alpha_2-\alpha_3+1}}\,dx\\
& = {e^{-4\pi}}.
\end{align*}
Finally, note that $L(s,\pi_{\mu^\vee},{\rm Ad}) = \Gamma_\R(s+1)^2\Gamma_\C(s+\alpha_2-\alpha_3)$.
Hence
\[
\frac{L(1,\pi_{\mu^\vee},{\rm Ad})}{\Gamma_\R(2)^2} = 2^{-\alpha_2+\alpha_3-2}(f_{\pi_{\mu^\vee}},f_{\pi_{\mu^\vee}}).
\]
This completes the proof.
\end{proof}

%\bibliographystyle{alpha}
%\bibliography{ref}

\begin{thebibliography}{AGI{\etalchar{+}}24}

\bibitem[AGI{\etalchar{+}}24]{AGIKMS2024}
H.~Atobe, W.~T. Gan, A.~Ichino, T.~Kaletha, A.~M\'{i}nguez, and S.~W. Shin.
\newblock {Local intertwining relations and co-tempered $A$-packets of
  classical groups}.
\newblock 2024.
\newblock arXiv:2410.13504.

\bibitem[AOY24]{AOY2024}
H.~Atobe, M.~Oi, and S.~Yasuda.
\newblock {Local newforms for generic representations of unramified odd unitary
  groups andfundamental lemma}.
\newblock {\em Duke Math. J.}, 173:2447--2479, 2024.

\bibitem[Ato20]{Atobe2020}
H.~Atobe.
\newblock {On the non-vanishing of theta liftings of tempered representations
  of ${\rm U}(p,q)$}.
\newblock {\em Adv. Math.}, 363, 2020.
\newblock 106984.

\bibitem[BG14]{BG2014}
K.~Buzzard and T.~Gee.
\newblock {The conjectural connections between automorphic representations and
  Galois representations}.
\newblock In {\em Automorphic Forms and Galois Representations, Volume 1},
  volume 414 of {\em London Mathematical Society Lecture Note Series}, pages
  135--187. Cambridge University Press, 2014.

\bibitem[BHR94]{BHR1994}
D.~Blasius, M.~Harris, and D.~Ramakrishnan.
\newblock {Coherent cohomology, limits of discrete seres, and Galois
  conjugation}.
\newblock {\em Duke Math. J.}, 73(3):647--685, 1994.

\bibitem[BPC25]{BPC2023}
R.~Beuzart-Plessis and P.-H. Chaudouard.
\newblock {The global Gan-Gross-Prasad conjecture for unitary groups. II. From
  Eisenstein series to Bessel periods}.
\newblock {\em Forum Math. Pi}, 13(16):1--98, 2025.

\bibitem[BR17]{BR2017}
B.~Balasubramanyam and A.~Raghuram.
\newblock {Special values of adjoint $L$-functions and congruences for
  automorphic forms on ${\rm GL}(n)$ over a number field}.
\newblock {\em Amer. J. Math.}, 139(3):641--679, 2017.

\bibitem[Che21]{Chen2021f}
S.-Y. Chen.
\newblock {On Deligne's conjecture for symmetric sixth $L$-functions of Hilbert
  modular forms}.
\newblock 2021.
\newblock arXiv:2110.06261.

\bibitem[Che22a]{Chen2021c}
S.-Y. Chen.
\newblock {Algebraicity of critical values of adjoint $L$-functions for ${\rm
  GSp_4}$}.
\newblock {\em Ramanujan J.}, 59:883--931, 2022.

\bibitem[Che22b]{Chen2020}
S.-Y. Chen.
\newblock {Algebraicity of the near central non-critical values of symmetric
  fourth $L$-functions for Hilbert modular forms}.
\newblock {\em J. Number Theory}, 231:269--315, 2022.

\bibitem[Che23a]{Chen2021}
S.-Y. Chen.
\newblock {On Deligne's conjecture for symmetric fourth $L$-functions of
  Hilbert modular forms}.
\newblock {\em Adv. Math.}, 414, 2023.
\newblock 108860.

\bibitem[Che23b]{Cheng2023}
Y.~Cheng.
\newblock {Local newforms for generic representations of unramified ${\rm
  U}_{2n+1}$ and Rankin-Selberg integrals}.
\newblock 2023.
\newblock arXiv:2207.02118.

\bibitem[CI23]{CI2019}
S.-Y. Chen and A.~Ichino.
\newblock {On Petersson norms of generic cusp forms and special values of
  adjoint $L$-functions for ${\rm GSp}_4$}.
\newblock {\em Amer. J. Math.}, 145:899--993, 2023.

\bibitem[Clo90]{Clozel1990}
L.~Clozel.
\newblock {Motifs et Formes Automorphes: Applications du Principe de
  Fonctorialit\'e}.
\newblock In {\em Automorphic Forms, Shimura Varieties, and L-functions, Vol.
  I}, Perspectives in Mathematics, pages 77--159, 1990.

\bibitem[CS80]{CS1980}
W.~Casselman and J.~A. Shalika.
\newblock {The unramified principal series of $p$-adic groups. II. The
  Whittaker function}.
\newblock {\em Compos. Math.}, 41:207--231, 1980.

\bibitem[Dim05]{Dimitrov2005}
M.~Dimitrov.
\newblock {Galois representations modulo $p$ and cohomology of Hilbert modular
  varieties}.
\newblock {\em Ann. Sci. \'Ec. Norm. Sup\'er.}, 38:505--551, 2005.

\bibitem[FM24]{FM2024}
M.~Furusawa and K.~Morimoto.
\newblock {On the Gross--Prasad conjecture with its refinement for $({\rm
  SO}(5),{\rm SO}(2))$ and the generalized Bocherer conjecture}.
\newblock {\em Compos. Math.}, 160:2115--2202, 2024.

\bibitem[Gha02]{Ghate2002}
E.~Ghate.
\newblock {Adjoint $L$-values and primes of congruence for Hilbert modular
  forms}.
\newblock {\em Compos. Math.}, 132:243--281, 2002.

\bibitem[GI16]{GI2016}
W.~T. Gan and A.~Ichino.
\newblock {The Gross-Prasad conjecture and local theta correspondence}.
\newblock {\em Invent. Math.}, 206:705--799, 2016.

\bibitem[GR13]{GR2013}
W.~T. Gan and A.~Raghuram.
\newblock {Arithmeticity for periods of automorphic forms}.
\newblock In {\em Automorphic Representations and $L$-functions}, pages
  187--229. Tata Institute of Fundamental Research, 2013.

\bibitem[Gro97]{Gross1997}
B.~Gross.
\newblock {On the motive of a reductive group}.
\newblock {\em Invent. Math.}, 130:287--313, 1997.

\bibitem[Har85]{Harris1985}
M.~Harris.
\newblock {Arithmetic vector bundles and automorphic forms on Shimura
  varieties. I}.
\newblock {\em Invent. Math.}, 82:151--189, 1985.

\bibitem[Har86]{Harris1986}
M.~Harris.
\newblock {Arithmetic vector bundles and automorphic forms on Shimura
  varieties, II}.
\newblock {\em Compos. Math.}, 60(3):323--378, 1986.

\bibitem[Har89]{Harris1989}
M.~Harris.
\newblock {Functorial properties of toroidal compactifications of locally
  symmetric varieties}.
\newblock {\em Proc. Lond. Math. Soc.}, 3(59):1--22, 1989.

\bibitem[Har90a]{Harris1990b}
M.~Harris.
\newblock Automorphic forms and the cohomology of vector bundles on shimura
  varieties.
\newblock In {\em Automorphic Forms, Shimura Varieties, and L-functions, Vol.
  II}, Perspectives in Mathematics, pages 41--91, 1990.

\bibitem[Har90b]{Harris1990}
M.~Harris.
\newblock {Automorphic forms of $\overline{\partial}$-cohomology type as
  coherent cohomology classes}.
\newblock {\em J. Differential Geom.}, 32:1--63, 1990.

\bibitem[Hid81]{Hida1981}
H.~Hida.
\newblock {Congruences of cusp forms and special values of their zeta
  functions}.
\newblock {\em Invent. Math.}, 63:225--261, 1981.

\bibitem[How89]{Howe1989}
R.~Howe.
\newblock {Transcending classical invariant theory}.
\newblock {\em J. Amer. Math. Soc.}, 2(3):535--552, 1989.

\bibitem[HP25]{HP2025}
A.~Horawa and K.~Prasanna.
\newblock {Motivic action for Siegel modular forms}.
\newblock {\em J. Eur. Math. Soc.}, 2025.
\newblock to appear.

\bibitem[HY23]{HY2023}
M.-L. Hsieh and S.~Yamana.
\newblock {Five-variable $p$-adic $L$-functions for $U(3) \times U(2)$}.
\newblock 2023.
\newblock arXiv:2311.17661.

\bibitem[Ich22]{Ichino2022}
A.~Ichino.
\newblock {Theta lifting for tempered representations of real unitary groups}.
\newblock {\em Adv. Math.}, 398, 2022.
\newblock DOI:10.1016/j.aim.2022.108188.

\bibitem[JPSS81]{JPSS1981}
H.~Jacquet, I.~I. Piatetski-Shapiro, and J.~Shalika.
\newblock {Conducteur des r\'epresentations du groupe lin\'eaire}.
\newblock {\em Math. Ann.}, 256:199--214, 1981.

\bibitem[JS81]{JS1981b}
H.~Jacquet and J.~A. Shalika.
\newblock {On Euler products and the classification of automorphic forms II}.
\newblock {\em Amer. J. Math.}, 103(4):777--815, 1981.

\bibitem[Kal13]{Kaletha2013}
T.~Kaletha.
\newblock {Genericity and contragredience in the local Langlands
  correspondence}.
\newblock {\em Algebra Number Theory}, 7(10):2447--2474, 2013.

\bibitem[Kna86]{Knapp1986}
A.~W. Knapp.
\newblock {\em {Representation Theory of Semisimple Groups. An Overview Based
  on Examples}}, volume~36 of {\em Princeton Mathematical Series}.
\newblock Princeton University Press, 1986.

\bibitem[Kud94]{Kudla1994}
S.~S. Kudla.
\newblock {Splitting metaplectic covers of dual reductive pairs}.
\newblock {\em Israel J. Math.}, 87:361--401, 1994.

\bibitem[Lan13]{Lan2013}
K.-W. Lan.
\newblock {\em {Arithmetic compactifications of PEL-type Shimura varieties}},
  volume~36 of {\em London Mathematical Society Monographs}.
\newblock Princeton University Press, 2013.

\bibitem[LM15]{LM2015}
E.~Lapid and Z.~Mao.
\newblock {A conjecture on Whittaker-Fourier coefficients of cusp forms}.
\newblock {\em J. Number Theory}, 146:448--505, 2015.

\bibitem[LM17]{LM2017}
E.~Lapid and Z.~Mao.
\newblock {On an analogue of the Ichino-Ikeda conjecture for Whittaker
  coefficients on the metaplectic group}.
\newblock {\em Algebra Number Theory}, 11(3):713--765, 2017.

\bibitem[LO23]{LO2018}
F.~Lemma and T.~Ochiai.
\newblock {Endoscopic congruences modulo adjoint $L$-values for ${\rm
  GSp}(4)$}.
\newblock {\em Kyoto J. Math.}, 63(2):281--333, 2023.

\bibitem[LS19]{LS2019}
J.-P. Labesse and J.~Schwermer.
\newblock {Central morphisms and cuspidal automorphicrepresentations}.
\newblock {\em J. Number Theory}, 205:170--193, 2019.

\bibitem[Mat13]{Matringe2013}
N.~Matringe.
\newblock {Essential Whittaker functions for GL(n)}.
\newblock {\em Doc. Math.}, 18:1191--1214, 2013.

\bibitem[Mil90]{Milne1990}
J.~S. Milne.
\newblock {Canonical models of (mixed) Shimura varieties and automorphic vector
  bundles}.
\newblock In {\em Automorphic Forms, Shimura Varieties, and L-functions, Vol.
  I}, Perspectives in Mathematics, pages 283--414, 1990.

\bibitem[Mok15]{Mok2015}
C.-P. Mok.
\newblock {\em {Endoscopic classification of representations of quasi-split
  unitary groups}}, volume 235 of {\em Memoirs of the American Mathematical
  Society}.
\newblock American Mathematical Society, 2015.

\bibitem[Mor22]{Morimoto2022}
K.~Morimoto.
\newblock {On a certain local identity for Lapid--Mao's conjecture and formal
  degree conjecture: even unitary group case}.
\newblock {\em J. Inst. Math. Jussieu}, 21(4):1107--1161, 2022.

\bibitem[Mor24]{Morimoto2024}
K.~Morimoto.
\newblock {Onn Ichino--Ikeda type formula of Whittaker periods for unitary
  groups}.
\newblock 2024.
\newblock arXiv:2403.19166.

\bibitem[Nam15]{Namikawa2015}
K.~Namikawa.
\newblock {On a congruence prime criterion for cusp forms on ${\rm GL}_2$ over
  number fields}.
\newblock {\em J. Reine Angew. Math.}, 707:149--207, 2015.

\bibitem[Pau00]{Paul2000}
A.~Paul.
\newblock {Howe correspondence for real unitary groups II}.
\newblock {\em Proc. Amer. Math. Soc.}, 128(10):3129--3136, 2000.

\bibitem[Pin89]{Pink1989}
R.~Pink.
\newblock {\em {Arithmetic compactification of mixed Shimura varieties}}.
\newblock PhD thesis, Rheinischen Friedrich-Wilhelms-Universit\"{a}t Bonn,
  1989.

\bibitem[Ree91]{Reeder1991}
M.~Reeder.
\newblock {Old forms on ${\rm GL}_n$}.
\newblock {\em Amer. J. Math.}, 113(5):911--930, 1991.

\bibitem[Ric25]{Ricoul2025}
T.~Ricoul.
\newblock {\em {Real quadratic base changes for ${\rm GL}_3$ and integral
  relations of automorphic periods}}.
\newblock PhD thesis, Université Sorbonne Paris Nord, 2025.

\bibitem[RS18]{RS2018}
A.~Raghuram and M.~Sarnobat.
\newblock Cohomological representations and functorial transfer from classical
  groups.
\newblock In {\em Cohomology of arithmetic groups}, volume 245 of {\em
  {Springer Proceedings in Mathematics and Statistics}}, pages 157--176.
  Springer, 2018.

\bibitem[Shi11]{Shin2011}
S.~W. Shin.
\newblock {Galois representations arising from some compact Shimura varieties}.
\newblock {\em Ann. of Math.}, 173:1645--1741, 2011.

\bibitem[Sou05]{Soudry2005}
D.~Soudry.
\newblock {On Langlands functoriality from classical groups to ${\rm GL}_n$}.
\newblock {\em Ast\'{e}risque}, 298:335--390, 2005.

\bibitem[Su24]{Su2018}
J.~Su.
\newblock {Coherent cohomology of Shimura varieties and automorphic forms}.
\newblock {\em Ann. Sci. \'Ec. Norm. Sup\'er.}, 57:1039--1099, 2024.

\bibitem[Urb95]{Urban1995}
E.~Urban.
\newblock {Formes automorphes cuspidales pour ${\rm GL}_2$ sur un corps
  quadratique imaginaire. Valeurs sp\'eciales de fonctions $L$ et congruences}.
\newblock {\em Compos. Math.}, 99(3):283--324, 1995.

\bibitem[Vog07]{Vogan2007}
D.~Vogan.
\newblock {Branching to a maximal compact subgroup}.
\newblock In {\em {Harmonic analysis, group representations, automorphic forms
  and invariant theory}}, volume~12 of {\em Lect. Notes Ser. Inst. Math. Sci.
  Natl.Univ. Singap.}, page 321–401, 2007.

\bibitem[Wal84]{Wallach1984}
N.~Wallach.
\newblock On the constant term of a square integrable automorphic form.
\newblock In {\em Operator algebras and group representations}, volume~18 of
  {\em {Monographs and Studies in Mathematics}}, pages 227--237, 1984.

\end{thebibliography}

\newcommand{\etalchar}[1]{$^{#1}$}

\end{document}